 \documentclass[preprint,12pt]{elsarticle}
\usepackage{amsmath}
\usepackage{graphicx}
\usepackage{algorithmic}
\usepackage{algorithm}
\usepackage{amssymb}
\usepackage{color}
\usepackage[tight,footnotesize]{subfigure}

\DeclareGraphicsExtensions{.eps}

\def\x{{\mathbf x}}
\def\v{{\mathbf v}}

\def\CF{{\mathcal{S}}}

\def\M {{\mathbf{M}}}         
\def\OM {{\mathbf{\Omega}}}    
\def\P {{\mathbf{P}}}         
\def\Q {{\mathbf{Q}}}         
\def\x {{\mathbf{x}}}        
\def\e {{\mathbf{e}}}        
\def\y {{\mathbf{y}}}        
\def\w {{\mathbf{w}}}        
\def\mdim {m}   
\def\sdim {d}   
\def\pdim {p}   
\def\cosp {\ell} 
\def \D {\mathbf{D}}         
\def \I {\mathbf{I}}         

\newtheorem{thm}{Theorem}[section]
\newtheorem{cor}[thm]{Corollary}
\newtheorem{lem}[thm]{Lemma}

\newtheorem{defn}[thm]{Definition}
\newtheorem{rem}[thm]{Remark}

\newcommand{\norm}[1]{\left\Vert#1\right\Vert}
\newcommand{\abs}[1]{\left\vert#1\right\vert}

\newcommand\vect[1]{{\bf#1}}
\newcommand\matr[1]{{\bf#1}}

\newcommand\alphabf{{\boldsymbol{\alpha}}}
\newcommand{\argmin}{\operatornamewithlimits{argmin}}
\newcommand{\argmax}{\operatornamewithlimits{argmax}}

\newcommand{\Real}{\mathbb R}
\newcommand\RR[1]{\mathbb{R}^{#1}}

\DeclareMathOperator{\supp}{supp}

\DeclareMathOperator{\range}{range}
\DeclareMathOperator{\rank}{rank}

\DeclareMathOperator{\conv}{conv}

\newcommand{\rg}[1]{\textcolor{black}{#1}}

\begin{document}

\begin{frontmatter}

\title{A Greedy Algorithm for the Analysis Transform Domain}
\author[csa]{R.~Giryes\corref{cor1}}

\cortext[cor1]{Corresponding author}

\address[csa]{The Department of Electrical and Computer Engineering,
        Duke University,
        Durham, NC 27708, USA}



\begin{abstract}
\rg{Many image processing applications benefited remarkably from the theory of sparsity.
One model of sparsity is the cosparse analysis one. It was shown that using $\ell_1$-minimization one might stably recover a cosparse signal from a small set of random linear measurements if the operator is a frame. Another effort has provided guarantees for dictionaries that have a near optimal projection procedure using greedy-like algorithms. However, no claims have been given for frames. A common drawback of all these existing techniques is their high computational cost for large dimensional problems.}

\rg{In this work we propose a new greedy-like technique with theoretical recovery guarantees for frames as the analysis operator, closing the gap between greedy and relaxation techniques.
Our results cover both the case of bounded adversarial noise, where we show that the algorithm provides us with a stable reconstruction, and the one of random Gaussian noise, for which we prove that it has a denoising effect, closing another gap in the analysis framework.
Our proposed program, unlike the previous greedy-like ones that solely act in the signal domain, operates mainly in the analysis operator's transform domain. 
Besides the theoretical benefit, the main advantage of this strategy is its computational efficiency that makes it easily applicable to visually big data. We demonstrate its performance on several high dimensional images.}

\end{abstract}


\end{frontmatter}



\section{Introduction}
\label{sec:intro}

For more than a decade the idea that signals may be represented sparsely
has a great impact on the field of signal and image processing.
New sampling theory has been developed \cite{Donoho06Stable} together with new tools for handling signals in different types of applications, such as image denoising \cite{dabov06image},
image deblurring \cite{Danielyan12BM3D}, super-resolution \cite{Zeyde12single}, radar \cite{Fannjiang10Compressed}, medical imaging \cite{Lustig08Compressed} and astronomy \cite{Salmon13Poisson}, to name a few \cite{Bruckstein09From}.
Remarkably, in most of these fields the sparsity based techniques achieve state-of-the-art results.

The classical sparse model is the synthesis one. In this model
the signal $\x \in \RR{d}$ is assumed to have a $k$-sparse representation $\alphabf \in \RR{n}$ under a given dictionary $\D \in \RR{d \times n}$. Formally,
\begin{eqnarray}
\x = \D \alphabf, ~~ \norm{\alphabf}_0 \le k,
\end{eqnarray}
where $\norm{\cdot}_0$ is the $\ell_0$-pseudo norm that counts the number of non-zero entries in a vector.
Notice, that the non-zero elements in $\alphabf$ corresponds to a set of columns that creates a low-dimensional subspace in which $\x$ resides.

Recently, a new sparsity based model has been introduced: the analysis one \cite{elad07Analysis, Nam12Cosparse}.
In this framework, we look at the coefficients of $\OM\x$, the coefficients of the signal after applying the transform $\OM \in \RR{p\times d}$ on it. The sparsity of the signal is measured by the number of zeros in $\OM\x$.
We say that a signal is $\cosp$-cosparse if  $\OM\x$ has $\cosp$ zero elements. Formally,
\begin{eqnarray}
\norm{\OM\x}_0 \le p-\cosp.
\end{eqnarray}

Remark that each zero element in $\OM\x$ corresponds to a row in $\OM$ to which the signal is orthogonal and all these rows define a subspace to which the signal is orthogonal. Similar to synthesis, when the number of zeros is large the signal's subspace is low dimensional.
Though the zeros are those that define the subspace, in some cases it is more convenient to use the number of non-zeros $k = p - \cosp$ as done in \cite{Candes11Compressed, Vaiter13Robust}.


The main setup in which the above models have been used is
\begin{eqnarray}
\label{eq:meas}
\y = \M \x + \e,
\end{eqnarray}
where $\y \in \RR{m}$ is a given set of measurements,  $\M \in \RR{m \times d}$ is the measurement matrix and $\e \in \RR{m}$ is an additive noise which is assumed to be either adversarial bounded noise \cite{Donoho06Stable,Bruckstein09From, Candes06Near,Candes05Decoding}
or with a certain given distribution such as Gaussian \cite{Candes07Dantzig,Bickel09Simultaneous,Giryes12RIP}.
The goal is to recover $\x$ from $\y$ and this is the focus of our work.
For details about other setups, the curious reader may refer to \cite{Harmany12Spiral,
Soltanolkotabi13Robust,Giryes14Sparsity, Giryes14SparsityBased, Pang13Ranking, Pang14Learning, Pang10Robust}.

Clearly, without a prior knowledge it is impossible to recover $\x$ from $\y$ in the case $m < d$, or have a significant denoising effect when $\e$ is random with a known distribution.
Hence, having a prior, such as the sparsity one, is vital for these tasks.
Both the synthesis and the analysis models lead to (different) minimization problems
that provide estimates for the original signal $\x$.

In synthesis, the signal is recovered by its representation, using
\begin{eqnarray}
\label{eq:l0_synthesis}
\hat\alphabf_{S-\ell_0} = \argmin_{\tilde\alphabf \in \RR{n}}\norm{\tilde\alphabf}_0 & s.t &
\norm{\y - \M\D\tilde{\alphabf}}_2 \le \lambda_{\e},
\end{eqnarray}
where $\lambda_{\e}$ is an upper bound for  $\norm{\e}_2$ if the noise is bounded and adversarial ($S-\ell_0$ refers to synthesis-$\ell_0$). Otherwise, it is a scalar dependent on the noise distribution \cite{Candes07Dantzig,Bickel09Simultaneous,Candes06Modern}.
The recovered signal is simply $\hat\x_{S-\ell_0} = \matr{D}\hat\alphabf_{S-\ell_0}$.
In analysis, we have the following minimization problem.
\begin{eqnarray}
\label{eq:l0_analysis}
\hat\x_{A-\ell_0} = \argmin_{\tilde\x \in \RR{d}}\norm{\OM\tilde\x}_0 & s.t &
\norm{\y - \M\tilde{\x}}_2 \le \lambda_{\e}.
\end{eqnarray}
The values of $\lambda_\e$ are selected as before depending on the noise properties ($A-\ell_0$ refers to analysis-$\ell_0$).


Both \eqref{eq:l0_synthesis} and \eqref{eq:l0_analysis} are NP-hard problems \cite{Nam12Cosparse, Davis97Adaptive}. Hence, approximation techniques are required. These are divided mainly into two categories: relaxation methods and
greedy algorithms.
In the first category we have the $\ell_1$-relaxation  \cite{elad07Analysis}  and the Dantzig selector \cite{Candes07Dantzig}, where the latter has been proposed  only for synthesis. The $\ell_1$-relaxation leads to the  following minimization problems
for synthesis and analysis respectively\footnote{\rg{Note that setting $\OM$ to be the finite difference operator in~\eqref{eq:l1_analysis} leads to the anisotropic version of the well-known total variation (TV) \cite{Rudin92TV}. See \cite{Needell13Stable,Giryes15Sampling} for more details.}}:
\begin{eqnarray}
\label{eq:l1_synthesis}
\hat\alphabf_{S-\ell_1} = \argmin_{\tilde\alphabf \in \RR{n}}\norm{\tilde\alphabf}_1 & s.t &
\norm{\y - \M\D\tilde{\alphabf}}_2 \le \lambda_{\e}, \\
\label{eq:l1_analysis}
\hat\x_{A-\ell_1} = \argmin_{\tilde\x \in \RR{d}}\norm{\OM\tilde\x}_1 & s.t &
\norm{\y - \M\tilde{\x}}_2 \le \lambda_{\e}.
\end{eqnarray}

Among the synthesis greedy strategies we mention the thresholding method,
orthogonal matching pursuit (OMP) \cite{Chen89Orthogonal,MallatZhang93, Pati93OMP}, CoSaMP \cite{Needell09CoSaMP}, subspace pursuit (SP) \cite{Dai09Subspace}, iterative hard thresholding \cite{Blumensath09Iterative} and hard thresholding pursuit (HTP) \cite{Foucart11Hard}. Their counterparts in analysis are thresholding \cite{Peleg12Performance}, GAP \cite{Nam12Cosparse}, analysis CoSaMP (ACoSaMP), analysis SP (ASP), analysis IHT (AIHT) and analysis HTP (AHTP) \cite{Giryes14Greedy}.

An important question to ask is what are the recovery guarantees that exist for these methods.
\rg{One main tool that was used for answering this question in the synthesis context is the restricted isometry property \cite{Candes06Near}.
It has been shown that under some conditions on 
the RIP of $\matr{MD}$, we have in the adversarial bounded noise case
 that
\begin{eqnarray}
\label{eq:rep_stable_bound}
\norm{\hat\alphabf_{alg} - \alphabf}_2^2 \le C_{alg}\norm{\e}_2^2,
\end{eqnarray} 
where $\hat\alphabf_{alg}$ is the recovered representation by one of the approximation algorithms and $C_{alg}>2$ is a constant that depends on the RIP of $\matr{MD}$ and differs for each of the methods \cite{Donoho06Stable,Candes06Near,
Needell09CoSaMP, Dai09Subspace, Blumensath09Iterative, Foucart11Hard, Zhang11Sparse, foucart10Sparse,Cai10New}. This result implies that these programs achieve a stable recovery.}

%

\rg{Similar results were provided for the case where the noise is random white Gaussian with variance $\sigma^2$ \cite{Candes07Dantzig,Bickel09Simultaneous,Giryes12RIP, BenHaim09Coherence}}. In this case the reconstruction error is guaranteed to be $O(k\log(n)\sigma^2)$ \cite{Candes07Dantzig, Bickel09Simultaneous, Giryes12RIP}. Unlike the adversarial noise case,
here we may have a denoising effect, as the recovery error can be smaller than the initial noise power $d\sigma^2$.
Remark that the above results can be extended also to the case where we have a model mismatch and the signal is not exactly $k$-sparse.

In the analysis framework we have similar guarantees for the adversarial noise case.
However, since the analysis model treats the signal directly, the guarantees are in terms of the signal and not the representation like in \eqref{eq:rep_stable_bound}. Two extensions for the RIP have been proposed providing guarantees for analysis algorithms.
The first is the D-RIP \cite{Candes11Compressed}:
\begin{defn}[D-RIP \cite{Candes11Compressed}]
\label{def:D_RIP}
A matrix $\matr{M}$ has the D-RIP with a dictionary $\matr{D}$ and a constant $\delta_k = \delta_{\matr{D},k}$,
if $\delta_{k}$ is the smallest constant that satisfies
\begin{eqnarray}
\label{eq:D_RIP}
&& \hspace{-0.5in} (1-\delta_{k})\norm{\matr{D}\tilde\alphabf}_2^2 \le \norm{\matr{M}\matr{D}\tilde\alphabf}_2^2 \le (1+\delta_{k})\norm{\matr{D}\tilde\alphabf}_2^2,
\end{eqnarray}
whenever $\tilde\alphabf$ is $k$-sparse.
\end{defn}
The D-RIP has been used for studying the performance of the analysis $\ell_1$-minimization \cite{Candes11Compressed, Liu12Compressed,Kabanava13Analysis}.
It has been shown that if $\OM$ is a frame with frame constants $A$ and $B$,
$\matr{D} = \OM^\dag$ and $\matr{M}$ has the D-RIP with  $\delta_{ak} \le \delta_{A-\ell_1}(a,A,B)$ then
\begin{eqnarray}
\label{eq:analysis_l1_stable_bound}
\norm{\hat\x_{A-\ell_1} - \x}_2^2 \le C_{A-\ell_1}\left(\norm{\e}_2^2 +
\frac{\norm{\OM\x - [\OM\x]_k}_1^2}{k}\right),
\end{eqnarray}
where the operator $[\cdot]_k$ is a hard thresholding operator that keeps the largest $k$ elements in a vector, $\delta_{A-\ell_1}(a,A,B)$ is a function of $a$, $A$ and $B$,
and $a \ge 1$ and $C_{A-\ell_1}$ are some constants.
A similar result has been proposed for analysis $\ell_1$-minimization with the finite difference operator \cite{Needell13Stable,Giryes15Sampling}.

\rg{The second is the O-RIP \cite{Giryes14Greedy}, which was used for the study of the greedy-like algorithms ACoSaMP, ASP, AIHT and AHTP.}
\begin{defn}[O-RIP \cite{Giryes14Greedy}]
\label{def:omega_RIP}
A matrix $\matr{M}$ has the O-RIP with an operator $\OM$ and a constant $\delta_{\OM,\cosp}$,
if $\delta_{\OM,\cosp}$ is the smallest constant that satisfies
\begin{eqnarray}
\label{eq:omega_RIP}
&& \hspace{-0.5in} (1-\delta_{\OM,\cosp})\norm{\vect{v}}_2^2 \le \norm{\matr{M}\vect{v}}_2^2 \le (1+\delta_{\OM,\cosp})\norm{\vect{v}}_2^2,
\end{eqnarray}
whenever $\OM\vect{v}$ has at least $\cosp$ zeroes.
\end{defn}
\rg{With the assumption} that there exists a cosupport selection procedure $\hat{\CF}_\cosp$ that implies a near optimal projection for $\OM$ with a constant $C_\ell$ (see Definition~\ref{def:C_optimal_proj} in Section~\ref{sec:analysis_alg}).
 It has  been proven for such operators that if $\delta_{\OM,a\ell}\le \delta_{alg}(C_\cosp, C_{2\cosp-p},\sigma_{\M}^2)$
then
\begin{eqnarray}
\label{eq:analysis_greedy_stable_bound}
\norm{\hat\x_{A-\ell_1} - \x}_2^2 \le C_{alg}\left(\norm{\e}_2^2 +
\norm{\x - \x^\cosp}_2^2
\right),
\end{eqnarray}
where $\sigma_{\M}^2$ is the largest singular value of $\M$, $\x^\cosp$ is the best $\cosp$-cosparse approximation for $\x$, $\delta_{alg}(C_\cosp, C_{2\cosp-p},\sigma_{\M}^2)$ is a function of $C_\cosp$, $C_{2\cosp-p}$ and $\sigma_{\M}^2$, and $a \ge 3$ and $C_{alg}$ are some constants that differ for each technique.

Notice that the conditions in synthesis imply that no linear dependencies can be allowed within small number of columns in the dictionary as the representation is the focus.
The existence of such dependencies may cause ambiguity in its recovery.
Since the analysis model performs in the signal domain, i.e. focus on the signal and not its representation, dependencies may be allowed within the dictionary. A recent series of contributions have shown that high correlations can be allowed in the dictionary also in the synthesis framework if the signal is the target and not the representation \cite{Blumensath11Sampling, Davenport13Signal, Giryes13CanP0, Giryes13OMP, Giryes13IHTconf, Giryes14GreedySignal, Giryes14NearOracle}.

\subsection{Our Contribution}

The conditions for greedy-like techniques require the constant $C_{\ell}$ to be close to $1$. Having a general projection scheme with $C_{\ell}=1$ is NP-hard \cite{Tillmann14Projection}.
The existence of a program with a constant close to one for a general operator is still an open problem. In particular, it is not known whether there exists a procedure that gives a small constant for frames. Thus, there is a gap between the results for the greedy techniques and the ones for the $\ell_1$-minimization.

Another drawback of the existing analysis greedy strategies is their high complexity. All of them require applying a projection to an analysis cosparse subspace, which implies a high computational cost.
Therefore, unlike in the synthesis case, they do not provide a ``cheap'' counterpart to the $\ell_1$-minimization.

In this work we propose a new efficient  greedy program, the transform domain IHT (TDIHT), which is an extension of IHT that operates in the analysis transform domain. Unlike AIHT, TDIHT has a low complexity, as it does not require applying computationally demanding projections like AIHT, and it inherits guarantees similar to
the ones of analysis $\ell_1$-minimization for frames.
Given that $\OM$ and its pseudo-inverse $\matr{D}$ can be applied efficiently, TDIHT demands in each iteration applying only $\OM$, $\matr{D}$, the measurement matrix $\matr{M}$ and its transpose together with other point-wise operations. This puts TDIHT as an efficient alternative to the existing analysis methods especially for big data problems, as for high dimensional problems all the existing techniques require solving computationally demanding high dimensional optimization problems \cite{Nam12Cosparse, Candes11Compressed, Giryes14Greedy}. The assumption about $\OM$ and $\matr{D}$ holds for many types of operators such as curvelet, wavelets, Gabor transforms and  discrete Fourier transform \cite{Candes11Compressed}.

Another gap exists between synthesis and analysis. To the best of our knowledge, no denoising guarantees has been proposed for analysis strategies  apart from the work in \cite{Peleg12Performance} that examines the performance of thresholding for the case $\M = \I$ and \cite{Giryes14effective} that studies the analysis $\ell_0$-minimization. We develop results for Gaussian noise in addition to the ones for adversarial noise, showing that it is possible to have a denoising effect using the analysis model also when $\M \ne \I$ and for different algorithms other than thresholding.

Our theoretical guarantees for TDIHT can be summarized by the following theorem:

\begin{thm}[Recovery Guarantees for TDIHT with Frames]
\label{thm:TDIHT_frame_guarantee}
Let $\y = \M\x + \e$ where $\norm{\OM\x}_0 \le k$ and $\OM$ is a  frame with frame constants $A$ and $B$, i.e.,  $\norm{\OM}_2 \le B$ and $\norm{\OM^\dag}_2 \le \frac{1}{A}$.
For certain selections of $\M$ and using only  $m = O(\frac{B}{A}k\log(p/k))$  measurements, the recovery result $\hat\x$ of TDIHT satisfies
\begin{eqnarray}
\label{eq:TDIHT_frame_guarantee_stablity}
\hspace{-0.1in} \norm{\x - \hat{\x}}_2 & \le &
 O\left(\frac{B}{A}\norm{\vect{e}}_2\right)
 \\ \nonumber && \hspace{-0.3in}
 + O\left(\frac{1+A}{A^2}\norm{\OM_{T^C}\x}_2 + \frac{1}{A^2\sqrt{k}}\norm{\OM_{T^C}\x}_1 \right),
\end{eqnarray}
for the case where $\e$ is an adversarial noise, implying that TDIHT leads to a stable recovery.
For the case that $\e$ is random i.i.d zero-mean Gaussian distributed noise with a known variance $\sigma^2$ we have
\begin{eqnarray}
\label{eq:TDIHT_frame_guarantee_denoising}
 \hspace{-0.1in} E\norm{\x - \hat{\x}}_2^2
&\le & O\left(\frac{B^2}{A^2}k\log(p)\sigma^2\right) \\ \nonumber && \hspace{-0.5in}
+  O\left(\frac{1+A}{A^2}\norm{\OM_{T^C}\x}_2 + \frac{1}{A^2\sqrt{k}}\norm{\OM_{T^C}\x}_1 \right)^2,
\end{eqnarray}
implying that TDIHT achieves a denoising effect.
\end{thm}

\begin{rem}
Note that $\OM_{T^C}\x = \OM\x - [\OM\x]_k$.
\end{rem}
\begin{rem}
Using Remark~2.3 in \cite{Needell09CoSaMP}, we can convert the $\ell_2$ norm into an $\ell_1$ norm in the model mismatch terms  in \eqref{eq:TDIHT_frame_guarantee_stablity} and \eqref{eq:TDIHT_frame_guarantee_denoising}, turning it to be more similar to what we have in the bound for analysis $\ell_1$-minimization in \eqref{eq:analysis_l1_stable_bound}.
\end{rem}
\begin{rem}
Theorem~\ref{thm:TDIHT_frame_guarantee} is a combination of  Theorems~\ref{thm:TDIHT_frame_adversarial_noise} and \ref{thm:TDIHT_frame_gaussian_noise}, plugging the minimal number of measurements implied by the D-RIP conditions of these theorems. Measurement matrices with sub-Gaussian entries are examples for matrices that satisfy this number of measurements \cite{Candes11Compressed}.
\end{rem}

\subsection{Organization}
This paper is organized as follows. Section~\ref{sec:notation} includes the notations used in this paper together with some preliminary lemmas for the D-RIP.  Section~\ref{sec:analysis_alg} presents the transform domain IHT. Section~\ref{sec:frame_guarantees} provides the proof of our main theorem. In Section~\ref{sec:exp} we present some simulations demonstrating the efficiency of TDIHT and its applicability to big data.
 Section~\ref{sec:conc}  concludes the paper.

\section{Notations and Preliminaries}
\label{sec:notation}

We use the following notation in our work:
\begin{itemize}
  \item We denote by $\norm{\cdot}_2$ the euclidean norm for vectors and the spectral ($2 \rightarrow 2$) norm for matrices; by $\norm{\cdot}_1$ the $\ell_1$ norm that sums the absolute values of a vector; and by $\norm{\cdot}_0$ the $\ell_0$ pseudo-norm which counts the number of non-zero elements in a vector.
  \item Given a cosupport set $\Lambda$, $\matr{\Omega}_{\Lambda}$ is a sub-matrix of $\matr{\Omega}$ with the {\em rows} that belong to $\Lambda$.

  \item In a similar way, for a support set $T$, $\matr{D}_T$ is a sub-matrix of $\matr{D}$
  with {\em columns}\footnote{By the abuse of notation we use the same notation for the selection sub-matrices of rows and columns. The selection will be clear from the context since in analysis the focus is always on the rows and in synthesis on the columns.} corresponding to the set of indices in $T$.
\item $(\w)_T$ keeps the elements in $\w$ supported on $T$ and zeros the rest.
\item $\supp(\cdot)$ returns the support of a vector and
        $\supp(\cdot, k)$ returns the support set of $k$ elements with largest magnitudes.
        \item $\left( \cdot \right)^\dag$ denotes the Moore-Penrose pseudo-inverse \cite{Moore1920PseudoInverse}.
  \item $\Q_\Lambda = \matr{I} - \OM_\Lambda^\dag\OM_\Lambda$ is the orthogonal projection onto the orthogonal
        complement of $\range(\matr{\Omega}_\Lambda^*)$.
  \item $\P_T = \matr{D}_T\matr{D}_T^\dag$ is the orthogonal projection onto $\range(\matr{D}_T)$.
\item Throughout the paper we assume that $n=p$.
  \item The original unknown $\cosp$-cosparse vector is denoted by $\x \in \RR{\sdim}$, its cosupport by $\Lambda$ and the support of the non-zero entries by $T = \Lambda^C$. By definition $| \Lambda | \ge \cosp$ and $| T | \le k$.
\item For a general $\cosp$-cosparse vector we use $\vect{v}\in \RR{d}$, for a general vector in the signal domain $\vect{z}\in \RR{d}$ and for a general vector in the analysis transform or dictionary representation domain $\w \in \RR{p}$.
\end{itemize}

We now turn to present several key properties of the D-RIP.
All of their proofs except of the last one, which we present hereafter,
appear in \cite{Giryes14GreedySignal}.

\begin{cor}
\label{cor:MP_RIP_norm}
If $\matr{M}$ satisfies the D-RIP with a constant $\delta_{k}$ then
\begin{eqnarray}
\label{eq:MP_RIP_norm}
\norm{\matr{M}\P_{T}}_2^2 \le 1+\delta_{k}
\end{eqnarray}
for every $T$ such that $\abs{T} \le k$.
\end{cor}

\begin{lem}
\label{lem:k_inequality}
For $k \le \tilde{k}$ it holds that $\delta_{k} \le \delta_{\tilde{k}}$.
\end{lem}

\begin{lem}
\label{lem:D_RIP_norm}
If $\matr{M}$ satisfies the D-RIP then
\begin{eqnarray}
\label{eq:D_RIP_norm}
\norm{\P_{T}(\matr{I} - \matr{M}^*\matr{M})\P_{T}}_2 \le \delta_k
\end{eqnarray}
for any $T$ such that $\abs{T} \le k$.
\end{lem}

\begin{cor}
\label{cor:D_RIP_norm_diff}
If $\matr{M}$ satisfies the D-RIP then
\begin{eqnarray}
\label{eq:D_RIP_norm_diff}
&& \norm{\P_{T_1}(\matr{I} - \matr{M}^*\matr{M})\P_{T_2}}_2 \le \delta_k, \\ \nonumber
\end{eqnarray}
for any $T_1$ and $T_2$ such that $\abs{T_1} \le k_1, \abs{T_2}\le k_2, k_1+k_2\le k$.
\end{cor}

The last Lemma we present is a generalization of Proposition~3.5 in \cite{Needell09CoSaMP}.
\begin{lem}
\label{lem:D_RIP_non_spars_up}
Suppose that $\M$ satisfies the upper inequality of the D-RIP, i.e.,
\begin{eqnarray}
\norm{\M\D\w}_2 \le \sqrt{1+\delta_k}\norm{\D\w}_2 && \forall \w, \norm{\w}_0 \le k,
\end{eqnarray}
and that $\norm{\matr{D}}_2 \le \frac{1}{A}$.
Then for any representation $\w$ we have
\begin{eqnarray}
\norm{\M\D\w}_2 \le \frac{\sqrt{1+\delta_k}}{A}\left(\norm{\w}_2 + \frac{1}{\sqrt{k}}\norm{\w}_1 \right).
\end{eqnarray}
\end{lem}
The proof is left to Appendix~\ref{sec:D_RIP_non_spars_up_proof}.
Before we proceed we recall the problem we aim at solving:
\begin{defn}[Problem $\cal{P}$]
Consider a measurement vector $\y \in \RR{\mdim}$
such that $\y=\M\x + \e$, where $\x\in \RR{\sdim}$ is either $\cosp$-cosparse
under a given and fixed analysis operator $\OM \in \RR{\pdim \times \sdim}$
or almost $\cosp$-cosparse, i.e. $\OM \x$ has $k = p - \cosp$ leading elements.
The non-zero locations of the $k$ leading elements is denoted by $T$.
$\M\in \RR{\mdim\times \sdim}$ is a degradation operator and $\e\in \RR{\mdim}$ is an additive noise.
Our task is to recover $\x$ from $\y$. The recovery result is denoted by $\hat\x$.
\end{defn}

\section{Transform Domain Iterative Hard Thresholding}
\label{sec:analysis_alg}

Our goal in this section is to provide a greedy-like approach that provide guarantees similar to the one of analysis $\ell_1$-minimization.
By analyzing the latter we notice that though it operates directly on the signal, it actually minimizes the coefficients in the transform domain. In fact, all the proof techniques utilized for this recovery strategy use the fact that nearness in the analysis dictionary domain implies the same in the signal domain \cite{Candes11Compressed, Needell13Stable, Liu12Compressed}.
Using this fact, recovery guarantees have been developed for tight frames \cite{Candes11Compressed}, general frames \cite{Liu12Compressed} and the 2D finite difference operator which corresponds to TV \cite{Needell13Stable}.
Working in the transform domain is not a new idea and was used before, especially in the context of dictionary learning \cite{Ophir11Multi,Ravishankar11MR,Ravishankar13Sparsifying}.

Henceforth, our strategy for extending the results of the $\ell_1$-relaxation is to modify the greedy-like approaches to operate in the transform domain.
In this paper we concentrate on IHT.
Before we turn to present the transform domain version of IHT we recall its synthesis and analysis versions.

\subsection{Quick Review of IHT and AIHT}
IHT and AIHT are assumed to know the cardinalities $k$ and $\cosp$ respectively. They aim at approximating variants of \eqref{eq:l0_synthesis}
and \eqref{eq:l0_analysis}:
\begin{eqnarray}
\label{eq:synthesisL0_k}
\argmin_{\tilde\alphabf} \norm{\y - \M \D \tilde\alphabf}_2^2 & s.t.& \norm{\tilde\alphabf}_0 \le k,
\end{eqnarray}
and
\begin{eqnarray}
\label{eq:analysisL0_k}
\argmin_{\tilde\x} \norm{\y - \M\tilde\x }_2^2 & s.t. & \norm{\OM\tilde\x}_0 \le p - \cosp.
\end{eqnarray}

IHT aims at recovering the representation and uses only one matrix $\matr{A} = \M\D$ in the whole recovery process. AIHT targets the signal and utilizes both $\M$ and $\OM$.  For recovering the signal using IHT, one has $\hat\x_{\text{\tiny IHT}} = \matr{D}\hat\alphabf_{\text{\tiny IHT}}$.
IHT \cite{Blumensath09Iterative} and AIHT \cite{Giryes14Greedy} are presented in Algorithms~\ref{alg:IHT} and \ref{alg:Analysis_IHT} respectively.

The iterations of IHT and AIHT are composed of two basic steps.
In both of them the first is a gradient step, with a step size $\mu_t$,
in the direction of minimizing $\norm{\y - \matr{M} \tilde\x}_2^2$.
The second step of IHT projects $\alphabf_g$ to the closest $k$-sparse subspace
by keeping the $k$ elements with the largest magnitudes. In AIHT a cosupport selection procedure, $\CF_\cosp$, is used for the cosupport selection and then orthogonal projection onto the corresponding orthogonal subspace is performed. In the theoretical study of AIHT this procedure is assumed to apply a near optimal projection:
\begin{defn}
\label{def:C_optimal_proj}
A procedure $\hat{\CF}_\cosp$ implies a near-optimal projection $\Q_{\hat{\CF}_\cosp(\cdot)}$ with a constant $C_\cosp$ if for any $\vect{z} \in \Real^d$
\begin{eqnarray}
\label{eq:C_optimal_proj}
&& \norm{\vect{z}-\Q_{\hat\CF_\cosp(\vect{z})}\vect{z}}_2^2 \le C_\cosp\min_{\Lambda}\norm{\vect{z} - \Q_{\Lambda}\vect{z}}_2^2.
\end{eqnarray}
\end{defn}
More details about this definition can be found in
\cite{Giryes14Greedy}.

In the algorithms' description we neither specify the stopping criterion, nor the step size selection technique. For exact details we refer the curious reader to
\cite{Giryes14Greedy, Blumensath09Iterative,Kyrillidis11Recipes,Cevher11ALPS}.

\begin{algorithm}
\caption{Iterative hard thresholding (IHT)} \label{alg:IHT}
\begin{algorithmic}

\REQUIRE $k, \matr{A}, \vect{y}$, where $\vect{y} = \matr{A}\alphabf
+ \vect{e}$, $k$ is the cardinality of $\alphabf$ and $\vect{e}$ is
an additive noise.

\ENSURE $\hat{\alphabf}_{\text{\tiny IHT}}$: $k$-sparse approximation of
$\alphabf$.

\STATE Initialize representation $\hat{\alphabf}^0 = \vect{0}$
 and set $t = 0$.

\WHILE{halting criterion is not satisfied}

\STATE $t = t + 1$.

\STATE Perform a gradient step: $\alphabf_g = \hat{\alphabf}^{t-1} + \mu^t \matr{A}^*(\y - \matr{A}\hat{\alphabf}^{t-1})$

\STATE Find a new support: $T^t = \supp(\alphabf_g,k)$

\STATE Calculate a new representation: $\hat{\alphabf}^t = (\alphabf_g)_{T^t}$.

\ENDWHILE

\STATE Form the final solution $\hat{\alphabf}_{\text{\tiny IHT}} = \hat{\alphabf}^t$.

\end{algorithmic}
\end{algorithm}

\begin{algorithm}
\caption{Analysis iterative hard thresholding (AIHT)} \label{alg:Analysis_IHT}
\begin{algorithmic}

\REQUIRE $\cosp, \matr{M}, \matr{\Omega}, \vect{y}$, where $\vect{y} = \matr{M}\vect{x}
+ \vect{e}$, $\cosp$ is the cosparsity of $\vect{x}$ under $\matr{\Omega}$ and $\vect{e}$ is
an additive noise.

\ENSURE $\hat{\vect{x}}_{\text{\tiny AIHT}}$: $\cosp$-cosparse approximation of
$\vect{x}$.

\STATE Initialize estimate $\hat{\x}^0 = \vect{0}$ and set $t = 0$.

\WHILE{halting criterion is not satisfied}

\STATE $t = t + 1$.

\STATE Perform a gradient step: $\vect{x}_g = \hat{\vect{x}}^{t-1} + \mu^t \M^*(\y - \M\hat{\x}^{t-1})$

\STATE Find a new cosupport: $\hat\Lambda^t = \hat{\CF}_\cosp(\x_g)$

\STATE Calculate a new estimate: $\hat{\x}^t = \Q_{\hat\Lambda^t}\vect{x}_g$.

\ENDWHILE

\STATE Form the final solution $\hat{\vect{x}}_{\text{\tiny AIHT}} = \hat{\vect{x}}^t$.

\end{algorithmic}
\end{algorithm}

\subsection{Transform Domain Analysis Greedy-Like Method}

The drawback of AIHT for handling analysis signals is that it assumes the existence of a near optimal cosupport selection scheme $\CF_\cosp$.
The type of analysis dictionaries for which a known feasible cosupport selection technique exists is very limited \cite{Giryes14Greedy,Tillmann14Projection}.
Note that this limit is not unique only to the analysis framework \cite{Blumensath11Sampling, Davenport13Signal,Giryes13IHTconf}.
Of course, it is possible to use a cosupport selection strategy with no guarantees on its near-optimality constant and it might work well in practice. For instance, simple hard thresholding has been shown to operate reasonably well in several instances where no practical projection is at hand \cite{Giryes14Greedy}. However, the theoretical performance guarantees
 depend heavily on the near-optimality constant.
 Since for many operators there are no known selection schemes with small constants, the existing guarantees for AIHT, as well as the ones of the other greedy-like algorithms,
 are very limited. In particular, to date, they do not cover frames and the 2D finite difference operator
 as the analysis dictionary.
Moreover, even when an efficient optimal selection procedure exists AIHT is required to apply a projection to the selected subspace which in many cases is computationally demanding.

To bypass these problems we propose an alternative greedy approach for the analysis framework that operates in the transform domain instead
of the signal domain. This strategy aims at finding the closest approximation to $\OM\vect{x}$ and not to $\vect{x}$ using the fact that for many analysis operators proximity in the transform domain implies the same in the signal domain.
In some sense, this approach imitates the classical synthesis techniques that recover the signal by putting the representation as the target.

In Algorithm~\ref{alg:TD_IHT} an extension for IHT for the transform domain is proposed.
This algorithm makes use of $k$, the number of non-zeros in $\OM \vect{x}$, and $\matr{D}$, a dictionary satisfying $\matr{D}\OM = \matr{I}$.
One option for $\matr{D}$ is $\matr{D} = \OM^\dag$. If $\OM$ does not have a full row rank, we may compute $\D$ by adding to $\OM$ rows that resides in its rows' null space and then applying the pseudo-inverse.
For example, for the 2D finite difference operator we may calculate $\D$ by computing its pseudo inverse with an additional row composed of ones divided by $\sqrt{d}$. However, this option is not likely to provide good guarantees.
As the 2D finite difference operator is beyond the scope of this paper we refer the reader to \cite{Needell13Stable} for more details on this subject.

Indeed, for many operators $\OM$, there are infinitely many options for $\D$ such that $\D\OM=\I$. When Algorithm~\ref{alg:TD_IHT} forms the final solution $\hat{\x}_{\text{\tiny TDIHT}}  = \D \hat{\w}^t$, different choices of $\D$ will lead to different results since $\hat{\w}^t$ may not belong to $\range(\OM)$. Therefore, in our study of the algorithm in the next section we focus on the choice of $\OM$ as a frame and $\D = \OM^\dag$.

For the step size $\mu_t$ we can use several options. The first one it use a constant step size $\mu_t = \mu$ for all iterations.
The second one is to use an `optimal' step size that minimizes the fidelity term
\begin{eqnarray}
\mu_t = \argmin_{\mu} \norm{\y - \M\D\hat{\w}^t}_2^2.
\end{eqnarray}
Since $\hat{\w}^t = (\w_g^t)_{\hat{T}^t}$ and $\hat{T}^t = \supp\left(\OM\matr{D}\hat{\w}^{t-1} + \mu^t \OM\M^*(\y - \M\matr{D}\hat{\w}^{t-1}),k\right)$ the above minimization problem does not have a close form solution and a line search for $\mu^t$ may involve a change of $\hat{T}^t$ several times. Instead, one may follow
\cite{Kyrillidis11Recipes} and limits the minimization to the support $\tilde{T} = \hat{T}^{t-1}  \cup  \supp\left( \OM\M^*(\y - \M\matr{D}\hat{\w}^{t-1}),k\right)$. In this case we have
\begin{eqnarray}
\hspace{-0.3in} \mu_t &=& \argmin_{\mu} \norm{\y - \M\matr{D}_{\tilde{T}}\OM_{\tilde{T}}\matr{D}\hat{\w}^{t-1} - \mu^t \M\D_{\tilde{T}}\OM_{\tilde{T}}\M^*(\y - \M\matr{D}\hat{\w}^{t-1})}_2^2 \\ \nonumber \hspace{-0.3in} &=& \frac{\left(\y - \M\matr{D}_{\tilde{T}}\OM_{\tilde{T}}\matr{D}\hat{\w}^{t-1}\right)^* \M\D_{\tilde{T}}\OM_{\tilde{T}}\M^*(\y - \M\matr{D}\hat{\w}^{t-1}) }{\norm{\M\D_{\tilde{T}}\OM_{\tilde{T}}\M^*(\y - \M\matr{D}\hat{\w}^{t-1})}_2^2}.
\end{eqnarray}
We denote this selection technique as the adaptive changing step-size selection.

In this paper we focus in the theoretical part on TDIHT with $\mu =1$.
Our results can be easily extended also to the other step-size selection options using the proof techniques in \cite{Giryes14Greedy,Kyrillidis11Recipes}.
Naturally, TDIHT with an adaptive changing step-size selection behaves better than TDIHT with a constant $\mu$. Therefore, in Section~\ref{sec:exp} we demonstrate the performance of the former.


\begin{algorithm}
\caption{Transform Domain Iterative hard thresholding (TDIHT)} \label{alg:TD_IHT}
\begin{algorithmic}

\REQUIRE $k, \M \in \RR{m \times d}, \OM \in \RR{p \times d}, \D \in \RR{d \times p}, \y$, where $\vect{y} = \matr{M}\vect{x}
+ \vect{e}$, $\D$ satisfies $\D\OM = \I$, $\x$ is $p-k$ cosparse under $\OM$ which implies that it has a $k$-sparse representation under $\D$, and $\e$ is
an additive noise.

\ENSURE $\hat{\vect{x}}_{\text{\tiny TDIHT}}$: Approximation of
$\vect{x}$.

\STATE Initialize estimate $\hat{\w}^0 = \vect{0}$ and set $t = 0$.

\WHILE{halting criterion is not satisfied}

\STATE $t = t + 1$.

\STATE Perform a gradient step: $\w_g^t = \OM\matr{D}\hat{\w}^{t-1} + \mu^t \OM\M^*(\y - \M\matr{D}\hat{\w}^{t-1})$

\STATE Find a new transform domain support: $\hat{T}^t = \supp(\w_g^t,k)$

\STATE Calculate a new estimate: $\hat{\w}^t = (\w_g^t)_{\hat{T}^t}$.

\ENDWHILE

\STATE Form the final solution $\hat{\vect{x}}_{\text{\tiny TDIHT}} = \matr{D}\hat{\w}^t$.

\end{algorithmic}
\end{algorithm}


\section{Frame Guarantees}
\label{sec:frame_guarantees}

We provide theoretical guarantees for the reconstruction performance of the transform domain analysis IHT (TDIHT), with a constant step size $\mu = 1$,  for  frames.
We start with the case that the noise is adversarial.

\begin{thm}[Stable Recovery of TDIHT with Frames]
\label{thm:TDIHT_frame_adversarial_noise}
Consider the problem $\cal P$ and apply TDIHT with a constant step-size $\mu = 1$ and $\D = \OM^\dag$.
Suppose that $\e$ is a bounded adversarial noise, $\OM$ is a frame with frame constants $0 < A, B < \infty$ such that $\norm{\OM}_2 \le B$ and $\norm{\matr{D}}_2 \le \frac{1}{A}$,
and $\matr{M}$ has the D-RIP with the dictionary $[\matr{D},\OM^*]$ and a constant $\delta_{\tilde{k}} = \delta_{[\matr{D},\OM^*],\tilde{k}}^D$.
If $\rho \triangleq \frac{2\delta_{ak}B}{A} < 1$ (i.e. $\delta_{ak} \le \frac{A}{2B}$), then after a finite number of iterations $t \ge t^* = \frac{\log\left((\sqrt{1+\delta_{2k}}B\norm{\vect{e}}_2+\eta)/\norm{\OM_T\x}_2\right)}{\log(\rho)}$
\begin{eqnarray}
\label{eq:TDIHT_frame_adversarial_noise}
 \hspace{-0.05in} \norm{\x - \hat{\x}^t}_2
&\le & \frac{2\eta}{A} + \frac{2B\sqrt{1+\delta_{2k}}}{A}\left(1+ \frac{1-\rho^t}{1-\rho} \right)\norm{\vect{e}}_2 \\ \nonumber && \hspace{-0.4in}  +
  \frac{2(1+\delta_{2k})}{A^2}\frac{1-\rho^t}{1-\rho}\bigg(\left(1+A+ \frac{A^2}{2}\right)\norm{\OM_{T^C}\x}_2  \\ \nonumber && ~~~~~~~~~~~~~~~~~~~~~~~~~~ + \frac{1}{\sqrt{k}}\norm{\OM_{T^C}\x}_1 \
   \bigg),
\end{eqnarray}
implying that TDIHT leads to a stable recovery.
For tight frames $a=3$ and for other frames $a=4$ .
\end{thm}

The result of this theorem is a generalization of the one presented in \cite{foucart10Sparse} for IHT.
Its proof follows from the following lemma.

\begin{lem}
\label{lem:TDIHT_xg_bound}
Consider the same setup of Theorem~\ref{thm:TDIHT_frame_adversarial_noise}.
Then the $t$-th iteration of TDIHT satisfies
\begin{eqnarray}
&& \hspace{-0.3in} \norm{\left(\OM\x - \w_g^t\right)_{T \cup \hat{T}^t }}_2 \le 2\delta_{ak}\frac{B}{A}\norm{\left(\OM\x - \w_g^{t-1}\right)_{T \cup \hat{T}^{t-1} }}_2
\\ \nonumber&& \hspace{-0.15in}  +\frac{{1+\delta_{2k}}}{A}\left(\left(1+A + \frac{A^2}{2}\right)\norm{\OM_{T^C}\x}_2
+ \frac{1}{\sqrt{k}}\norm{\OM_{T^C}\x}_1 \right)
\\ \nonumber&&  \hspace{-0.15in}  + \norm{\OM_{T \cup  \hat{T}^t}\matr{M}^*\vect{e}}_2.
\end{eqnarray}
\end{lem}

The proof of the above lemma is left to
Appendix~\ref{sec:TDIHT_xg_bound_proof}.
We turn now to prove Theorem~\ref{thm:TDIHT_frame_adversarial_noise}.

{\em Proof:}[Proof of Theorem~\ref{thm:TDIHT_frame_adversarial_noise}]
First notice that $\hat{\w}^0=0$ implies that $\w_g^0 =0$.
Using Lemma~\ref{lem:TDIHT_xg_bound}, recursion and the definitions of $\rho$ and $T_{\e} \triangleq \argmax_{\tilde{T}: \abs{\tilde{T}}\le k} \norm{\OM_{T\cup \tilde{T}}\M^*\e}_2$, we have that after $t$ iterations
\begin{eqnarray}
\label{eq:OMx_zgt_1st_bound}
&& \hspace{-0.3in} \norm{\left(\OM\x - \w_g^t\right)_{T \cup \hat{T}^t }}_2 \le \rho^t\norm{\left(\OM\x - \w_g^{0}\right)_{T \cup \hat{T}^{0} }}_2
\\ \nonumber && \hspace{-0.2in}
+  (1+\rho+\rho^2+\dots + \rho^{t-1})
  \bigg(  \norm{\OM_{T \cup  T_{\e}}\matr{M}^*\vect{e}}_2  \\ \nonumber &&
 ~~~~ + \frac{{1+\delta_{2k}}}{A}\left((1+A)\norm{\OM_{T^C}\x}_2 + \frac{1}{\sqrt{k}}\norm{\OM_{T^C}\x}_1 \right)
   \bigg)
  \\ \nonumber && \hspace{-0.3in}= \rho^t\norm{\OM_T\x}_2
 +
  \frac{1-\rho^t}{1-\rho}\bigg( \norm{\OM_{T \cup  T_{\e}}\matr{M}^*\vect{e}}_2 +   \\ \nonumber && ~~~~~~~
  \frac{{1+\delta_{2k}}}{A}\left((1+A)\norm{\OM_{T^C}\x}_2 + \frac{1}{\sqrt{k}}\norm{\OM_{T^C}\x}_1 \right)
   \bigg),
\end{eqnarray}
where the last equality is due to the equation of geometric series ($\rho <1$) and the facts that $\w_g^0 =0$ and $T^0 = \emptyset$.
For a given precision factor $\eta$ we have that if $t \ge t^* = \frac{\log\left((\norm{\OM_{T \cup  T_{\e}}\matr{M}^*\vect{e}}_2+\eta)/\norm{\OM_T\x}_2\right)}{\log(\rho)}$ then
\begin{eqnarray}
\label{eq:rhot_normOMTx_ineq}
\rho^t \norm{\OM_T\x}_2 \le \eta + \norm{\OM_{T \cup  T_{\e}}\matr{M}^*\vect{e}}_2.
\end{eqnarray}
As $\x = \D\OM\x$ and $\hat{\x}^t = \D \hat{\w}^t$, we have using matrix norm inequality
\begin{eqnarray}
\label{eq:D_OM_Frame_OMe_bound}
&& \hspace{-0.4in} \norm{\x - \hat{\x}^t}_2 = \norm{\D(\OM\x - \hat{\w}^t)}_2 \le
\frac{1}{A}\norm{\OM\x - \hat{\w}^t}_2.
\end{eqnarray}
Using the triangle inequality and the facts that $\hat{\w}^t$ is supported on $\hat{T}^t$ and $\norm{\OM_{(T \cup \hat{T}^t)^C}\x}_2 \le \norm{\OM_{T^C}\x}_2$, we have
\begin{eqnarray}
\norm{\x - \hat{\x}^t}_2 \le \frac{1}{A}\norm{\OM_{T^C}\x}_2
+\frac{1}{A}\norm{\left(\OM\x - \hat{\w}^t\right)_{T \cup \hat{T}^t}}_2.
\end{eqnarray}
By using again the triangle inequality and the fact that $\hat{\w}^t$ is the best $k$-term
approximation for $\w_g^t$ we get
\begin{eqnarray}
\label{eq:x_xhat_OMxWg_ineq}
\hspace{-0.2in} \norm{\x - \hat{\x}^t}_2 &\le & \frac{1}{A}\norm{\OM_{T^C}\x}_2 +
\frac{1}{A}\norm{(\OM\x - \w_g^t)_{T \cup \hat{T}^t}}_2
\\ \nonumber && + \frac{1}{A}\norm{(\w_g^t - \hat{\w}^t)_{T \cup \hat{T}^t}}_2
\\ \nonumber & \le &
 \frac{1}{A}\norm{\OM_{T^C}\x}_2 +
\frac{2}{A}\norm{(\OM\x - \w_g^t)_{T \cup \hat{T}^t}}_2.
\end{eqnarray}
Plugging \eqref{eq:rhot_normOMTx_ineq} and \eqref{eq:OMx_zgt_1st_bound} in \eqref{eq:x_xhat_OMxWg_ineq} yields
\begin{eqnarray}
\label{eq:x_xt_OMe_bound}
&& \hspace{-0.3in} \norm{\x - \hat{\x}^t}_2
\le \frac{2\eta}{A} + \frac{2}{A}\left(1+ \frac{1-\rho^t}{1-\rho} \right)\norm{\OM_{T \cup  T_{\e}}\matr{M}^*\vect{e}}_2 \\ \nonumber && +
  \frac{2(1+\delta_{2k})}{A^2}\frac{1-\rho^t}{1-\rho}\bigg(\left(1+A + \frac{A^2}{2}\right)\norm{\OM_{T^C}\x}_2 \\ \nonumber && ~~~~~~~~~~~~~~~~~~~~~~~~~~~~~~~~~~~~~~ + \frac{1}{\sqrt{k}}\norm{\OM_{T^C}\x}_1 \bigg).
\end{eqnarray}
Using the D-RIP and the fact that $\OM$ is a frame we have that
$\norm{\OM_{T \cup  T_{\e}}\matr{M}^*\vect{e}}_2 \le B\sqrt{1+\delta_{2k}}\norm{\vect{e}}_2$ and this completes the proof.
\hfill $\Box$
\bigskip

Having the result for the adversarial noise case we turn to give a bound for the case where a distribution of the noise is given.
We dwell on the white Gaussian noise case.
For this type of noise we make use of the following lemma.

\begin{lem}
\label{lem:gauss_noise_bound}
If $\vect{e}$ is a zero-mean white Gaussian noise with a variance $\sigma^2$ then
\begin{eqnarray}
E[\max_{\abs{T}\le k}\norm{\OM_{T}\M^*\e}_2^2] \le 4\max_i{\norm{\OM_i^*}}_2^2(1+\delta_1)k\log(p)\sigma^2.
\end{eqnarray}
\end{lem}
{\em Proof:}
First notice that the $i$-th entry in $\OM_{T}\M^*\e$ is Gaussian distributed random variable with zero-mean and variance $\norm{\M\OM_i^*}_2^2\sigma^2$.
Denoting by $\matr{W}$ a diagonal matrix such that $\matr{W}_{i,i} = \norm{\M\OM_i^*}_2$, we have that each entry in $\matr{W}^{-1}\OM_{T}\M^*\e$ is Gaussian distributed with variance $\sigma^2$. Therefore, using Theorem~2.4 from \cite{Giryes12RIP} we have
$E[\max_{\abs{T}\le k}\norm{\matr{W}^{-1}\OM_{T}\M^*\e}_2^2] \le 4k\log(p)\sigma^2$.
Using the D-RIP we have that $\norm{\M\OM_i^*}_2^2 \le (1+\delta_1)\norm{\OM_i^*}_2^2$. Since the $\ell_2$ matrix norm
of a diagonal matrix is the maximal diagonal element we have that $\norm{\matr{W}}_2^2 \le \max_i{(1+\delta_1)\norm{\OM_i^*}_2^2}$
and this provides the desired result.
\hfill $\Box$
\bigskip

\begin{thm}[Denoising Performance of TDIHT with Frames]
\label{thm:TDIHT_frame_gaussian_noise}
Consider the problem $\cal P$ and apply TDIHT with a constant step-size $\mu = 1$.
Suppose that $\e$ is an additive white Gaussian noise with a known variance $\sigma^2$ (i.e. for each element $\e_i \sim N(0,\sigma^2)$,
$\OM$ is a frame with frame constants $0 < A, B < \infty$ such that $\norm{\OM}_2 \le B$ and $\norm{\matr{D}}_2 \le \frac{1}{A}$,
and $\matr{M}$ has the D-RIP with the dictionary $[\matr{D},\OM^*]$ and a constant $\delta_{\tilde{k}} = \delta_{[\matr{D},\OM^*],\tilde{k}}^D$.
If $\rho \triangleq \frac{2\delta_{ak}B}{A} < 1$ (i.e. $\delta_{ak} \le \frac{A}{2B}$), then after a finite number of iterations $t \ge t^* = \frac{\log\left((\norm{\OM_{T \cup  T_{\e}}\matr{M}^*\vect{e}}_2+\eta)/\norm{\OM_T\x}_2\right)}{\log(\rho)}$
\begin{eqnarray}
 && \hspace{-0.31in} E\norm{\x - \hat{\x}^t}_2^2
\le  \frac{32(1+\delta_1)B^2}{A^2}\left(1+ \frac{1-\rho^t}{1-\rho} \right)^2k\log(p)\sigma^2 \\ \nonumber &&+\frac{8(1+\delta_{2k})^2}{A^4}\left(\frac{1-\rho^t}{1-\rho}\right)^2\bigg(\left(1+A + \frac{A^2}{2}\right)\norm{\OM_{T^C}\x}_2
\\ \nonumber &&
~~~~~~~~~~~~~~~~~~~~~~~~~~~~~~~~~~~~~~~~~ + \frac{1}{\sqrt{k}}\norm{\OM_{T^C}\x}_1 \bigg)^2,
\end{eqnarray}
implying that TDIHT has a denoising effect.
For tight frames $a=3$ and for other frames $a=4$.
\end{thm}

{\em Proof:}
Using the fact that for tight frames $\max_i{\norm{\OM_i^*}}_2^2 \le B$,
we have using Lemma~\ref{lem:gauss_noise_bound}
that
\begin{eqnarray}
E\norm{\OM_{T \cup  T_{\e}}\matr{M}^*\vect{e}}_2^2 \le
4B^2(1+\delta_1)k\log(p)\sigma^2.
\end{eqnarray}
Plugging this in a squared version of \eqref{eq:x_xt_OMe_bound} (with $\eta = 0$), using the fact that for any two constants $a,b$ we have $(a+b)^2 \le 2a^2 +2b^2$, leads to the desired result.
\hfill $\Box$
\bigskip


\subsection{Result's Optimality}
\label{sec:opt}

Notice that in all our conditions, the number of measurements we need is
$O((p-\cosp)\log(p))$ and it is not dependent explicitly on the intrinsic dimension of the signal.
Intuitively, we would expect the minimal number of measurements to be rather a function of $d-r$, where $r= \rank(\OM_\Lambda)$ (refer to \cite{Nam12Cosparse, Giryes14Greedy} for more details),
and henceforth our recovery conditions seems to be sub-optimal.

Indeed, this is the case with the analysis $\ell_0$-minimization problem \cite{Nam12Cosparse}. However, all the guarantees developed for feasible programs \cite{Candes11Compressed, Giryes14Greedy, Needell13Stable,  Liu12Compressed, Kabanava13Analysis}
require at least $O((p-\cosp)\log(p/(p-\cosp)))$ measurements.
Apparently, such conditions are too demanding because $d-r$, which can be much smaller than $p-\cosp$, does not play any role in them.
However, it seems that it is impossible to robustly reconstruct the signal with fewer measurements \cite{Giryes14effective}.

The same argument may be used for the denoising bound we have for TDIHT ($O((p-\cosp)\log(p)\sigma^2)$), saying that we would expect it to be $O((d - r)\log(p)\sigma^2)$. Interestingly, even the $\ell_0$ solution cannot achieve the latter bound but only a one which is a function of the cosparsity \cite{Giryes14effective}.

In this work we developed recovery guarantees for TDIHT with frames.
These guarantees close a gap between relaxation-based techniques and greedy algorithms in the analysis framework, and extend the denoising guarantees of synthesis methods to analysis.
It is interesting to ask whether these can be extended to other operators such as the 2D finite difference or to other methods such as AIHT or an analysis version of the Dantzig selector. The core idea in this work is the connection between the signal domain and the transform domain.
We believe that the relationships used in this work can be developed further, leading to other new results and improving existing techniques.
Some results in this direction appear in \cite{Giryes15Sampling}.

\section{Experiments}
\label{sec:exp}

In this section we repeat some of the experiments performed in \cite{Nam12Cosparse,Giryes14Greedy,Nam11GAPN}.
We start with synthetic signals in the noiseless case. We test the performance of TDIHT with an adaptive changing step-size and compare it with AIHT, AHTP, ASP, ACoSaMP \cite{Giryes14Greedy}, $\ell_1$-minimization \cite{elad07Analysis} and GAP \cite{Nam12Cosparse}.
As there are several possibilities for setting the parameters in AIHT, AHTP, ASP and ACoSaMP we select the ones that provide the best performance according to \cite{Giryes14Greedy}.

We draw a phase transition diagram \cite{Donoho09countingfaces} for each of the algorithms.
We use a random matrix $\matr{M}$ and a random tight frame $\OM$ with $d=120$ and $p=144$,
where each entry in the matrices is drawn independently from the Gaussian distribution.
We test $20$ different possible values of $m$ and $20$ different values of $\cosp$ and for each pair repeat the experiment $50$ times.
In each experiment we check whether we have a perfect reconstruction. White cells in the diagram
denote a perfect reconstruction in all the experiments of the pair and black cells denotes total failure in the reconstruction.
The values of $m$ and $\cosp$ are selected according to the following formula:
\begin{eqnarray}
m = \delta d && \cosp = d - \rho m,
\end{eqnarray}
where $\delta$, the sampling rate, is the x-axis of the phase diagram and
$\rho$, the ratio between the cosparsity of the signal and the number of measurements, is the y-axis.

\begin{figure*}[!t]
{\subfigure[AIHT, adaptive step-size]{\includegraphics[width=1.6in]{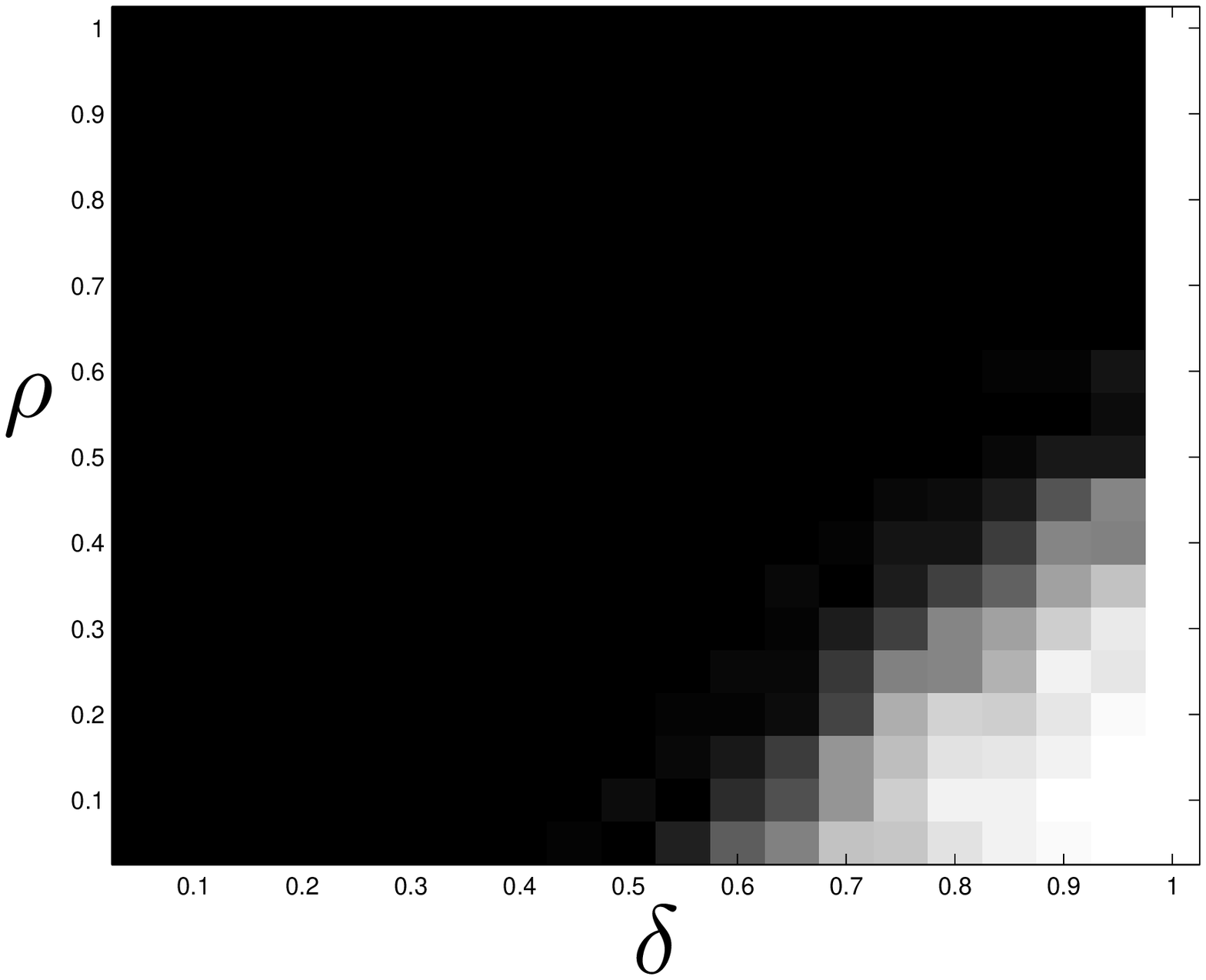}}}%
\hfil
\subfigure[AHTP, adaptive step-size]{\includegraphics[width=1.6in]{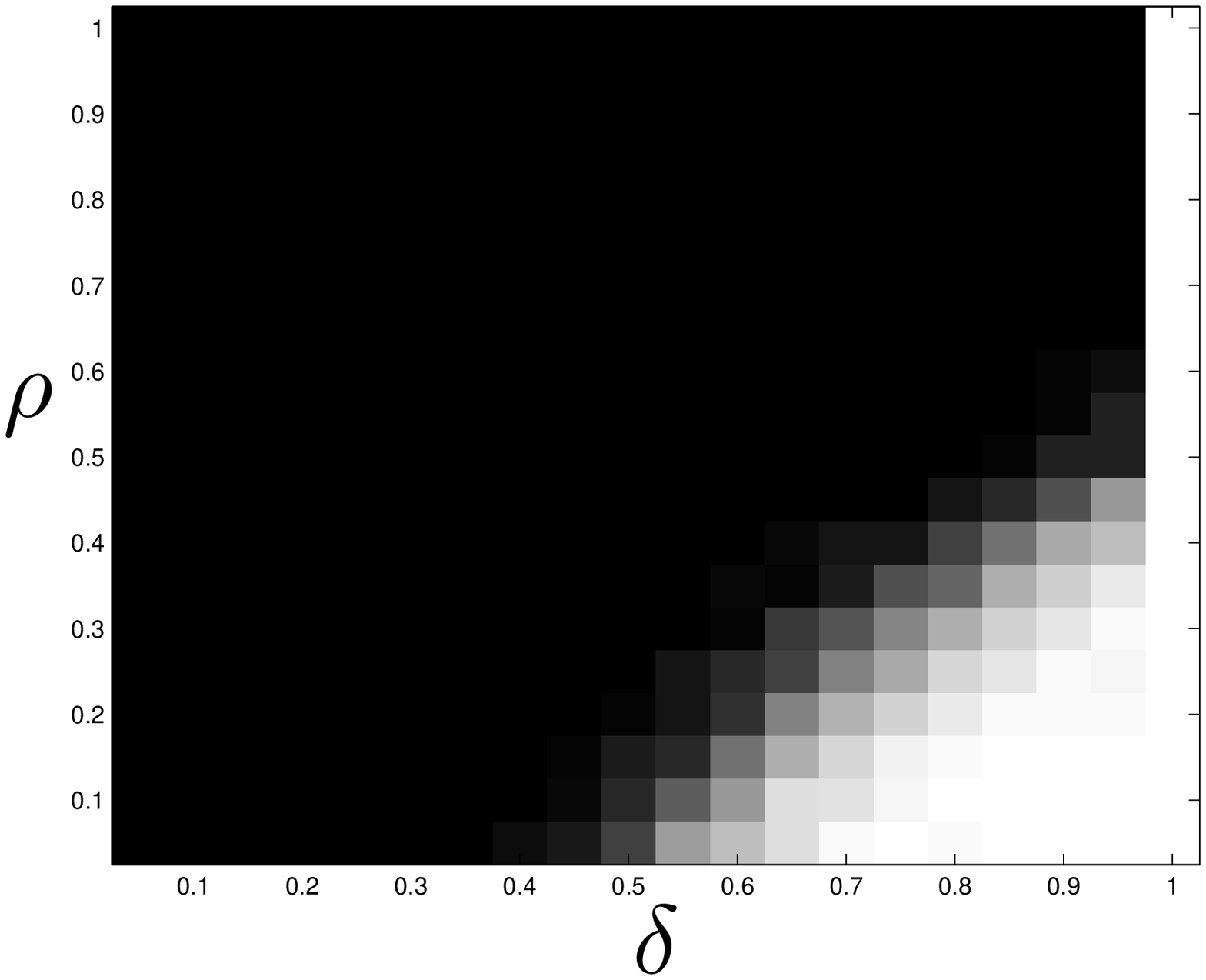}}%
\hfil
\subfigure[ACoSaMP]{\includegraphics[width=1.6in]{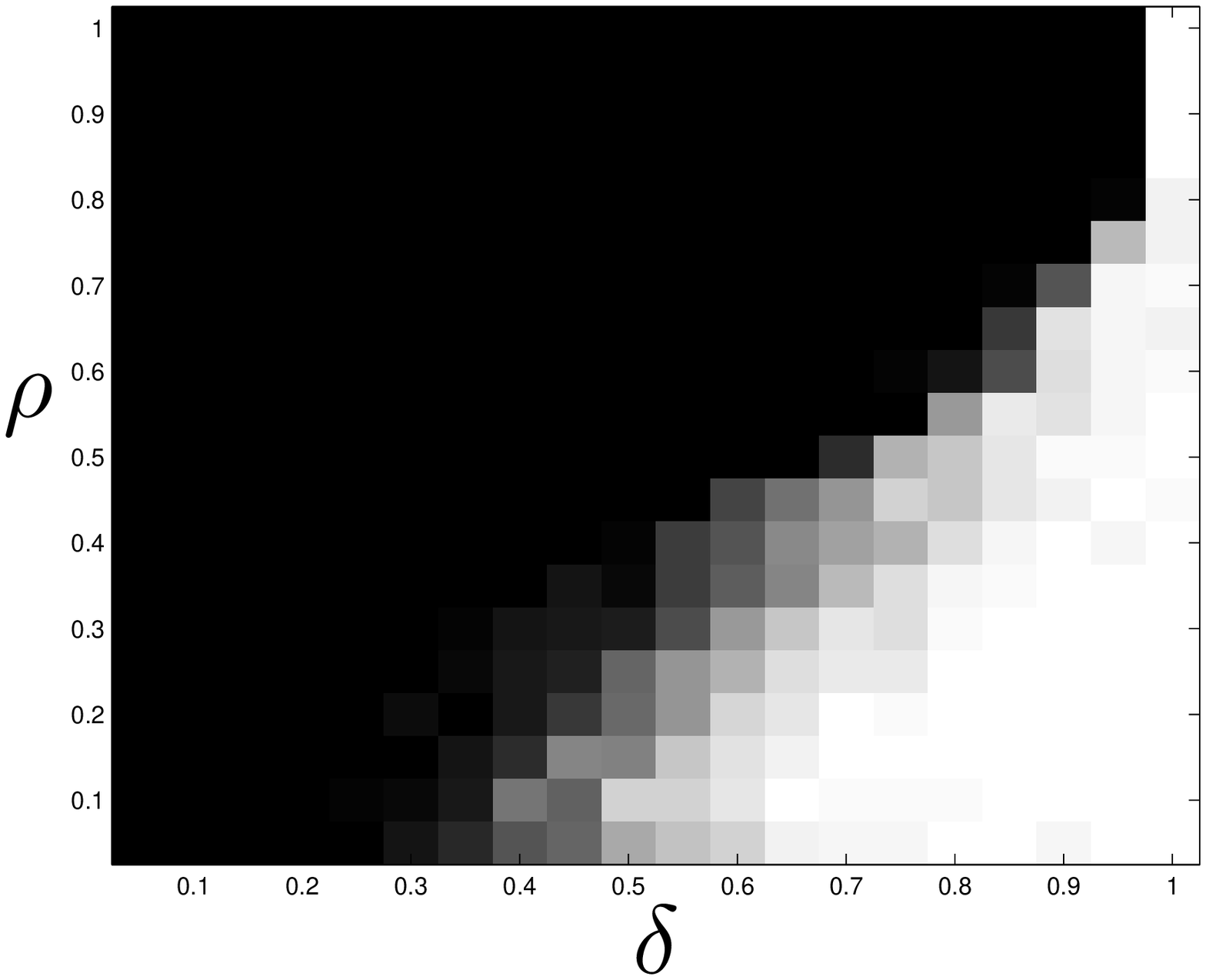}}%
\hfil
\subfigure[ASP]{\includegraphics[width=1.6in]{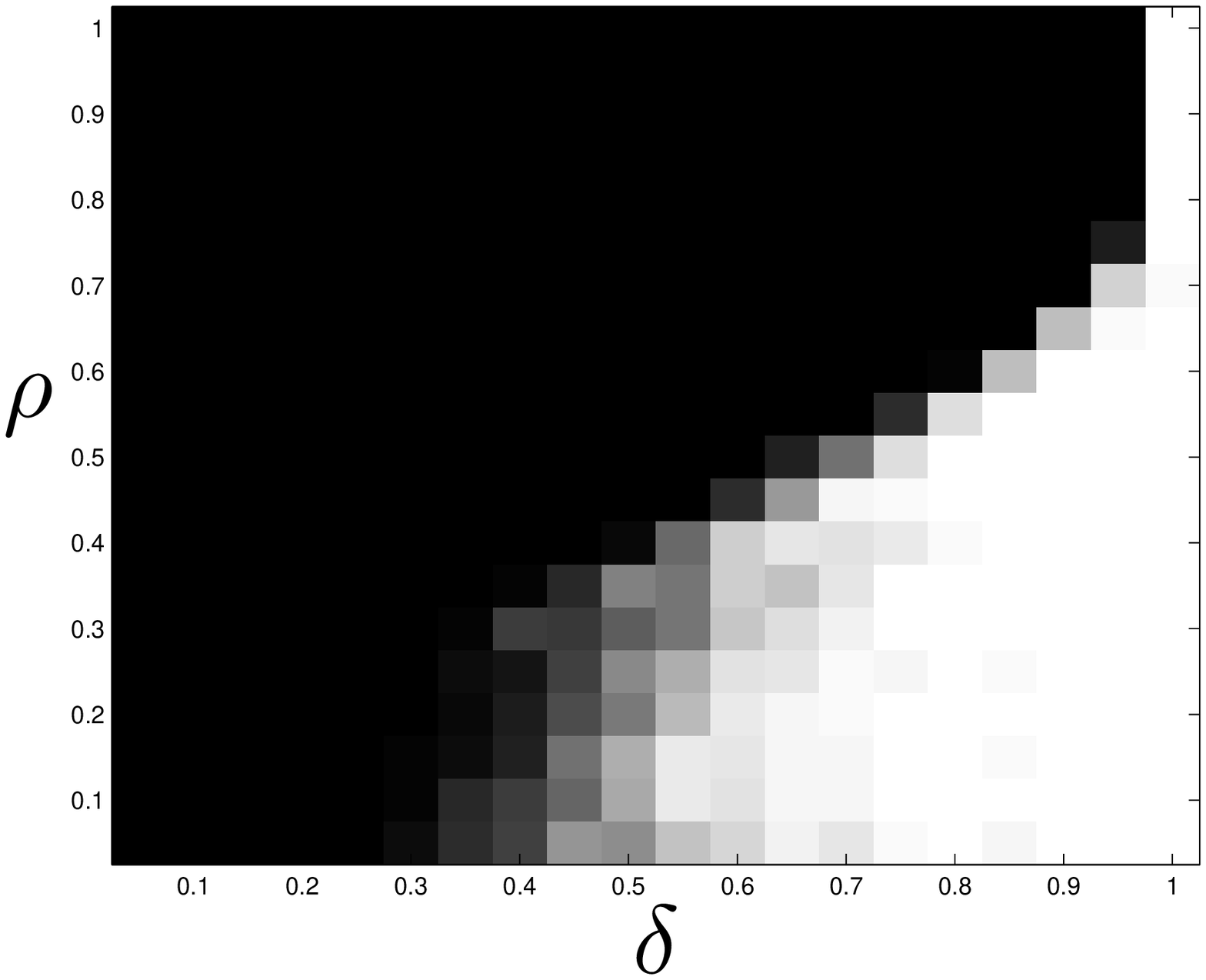}}%
\vfil
\subfigure[TDIHT, adaptive step-size]{\includegraphics[width=1.6in]{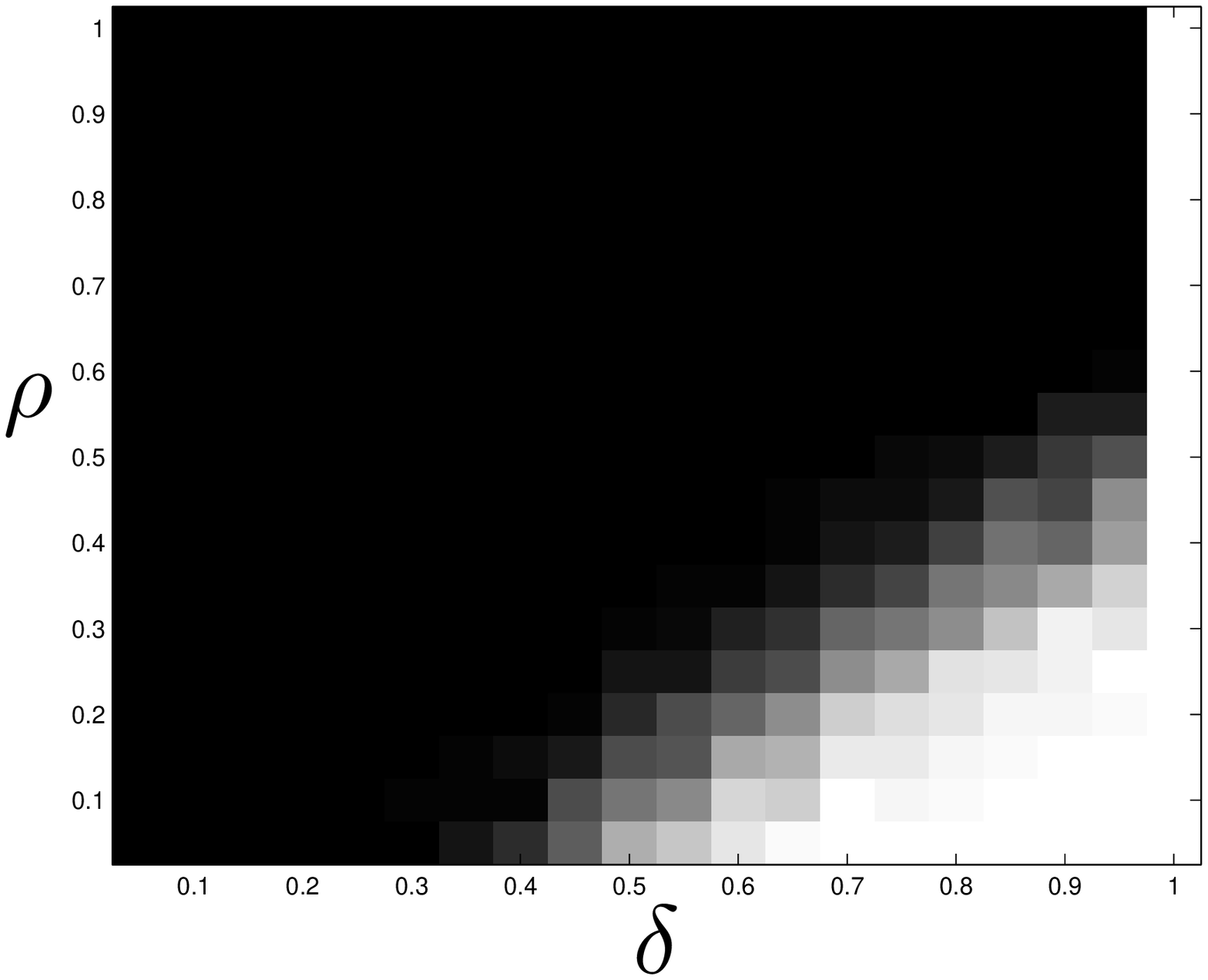}}%
\hfil
\subfigure[A-$\ell_1$-minimization]{\includegraphics[width=1.6in]{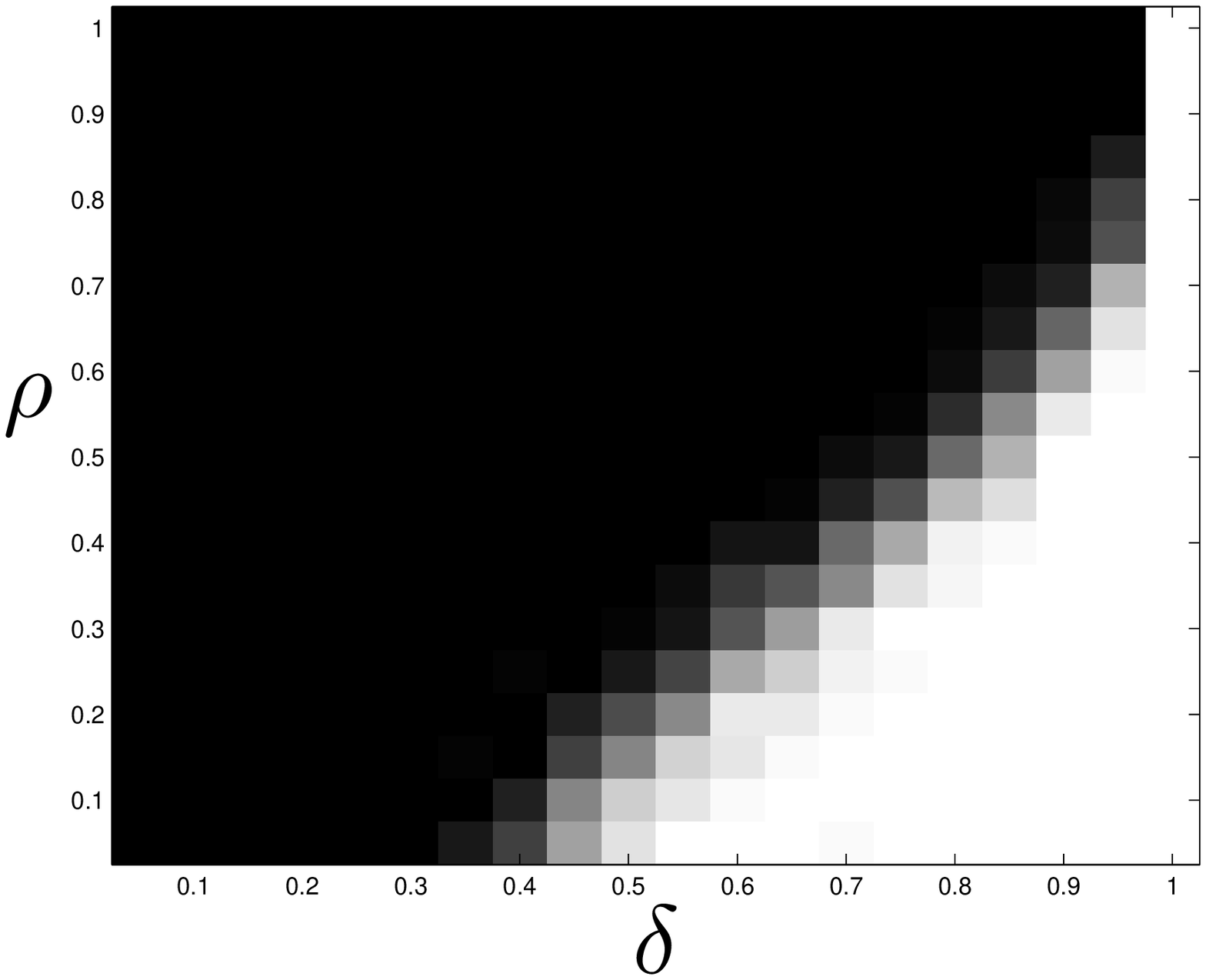}}%
\hfil
\subfigure[GAP]{\includegraphics[width=1.6in]{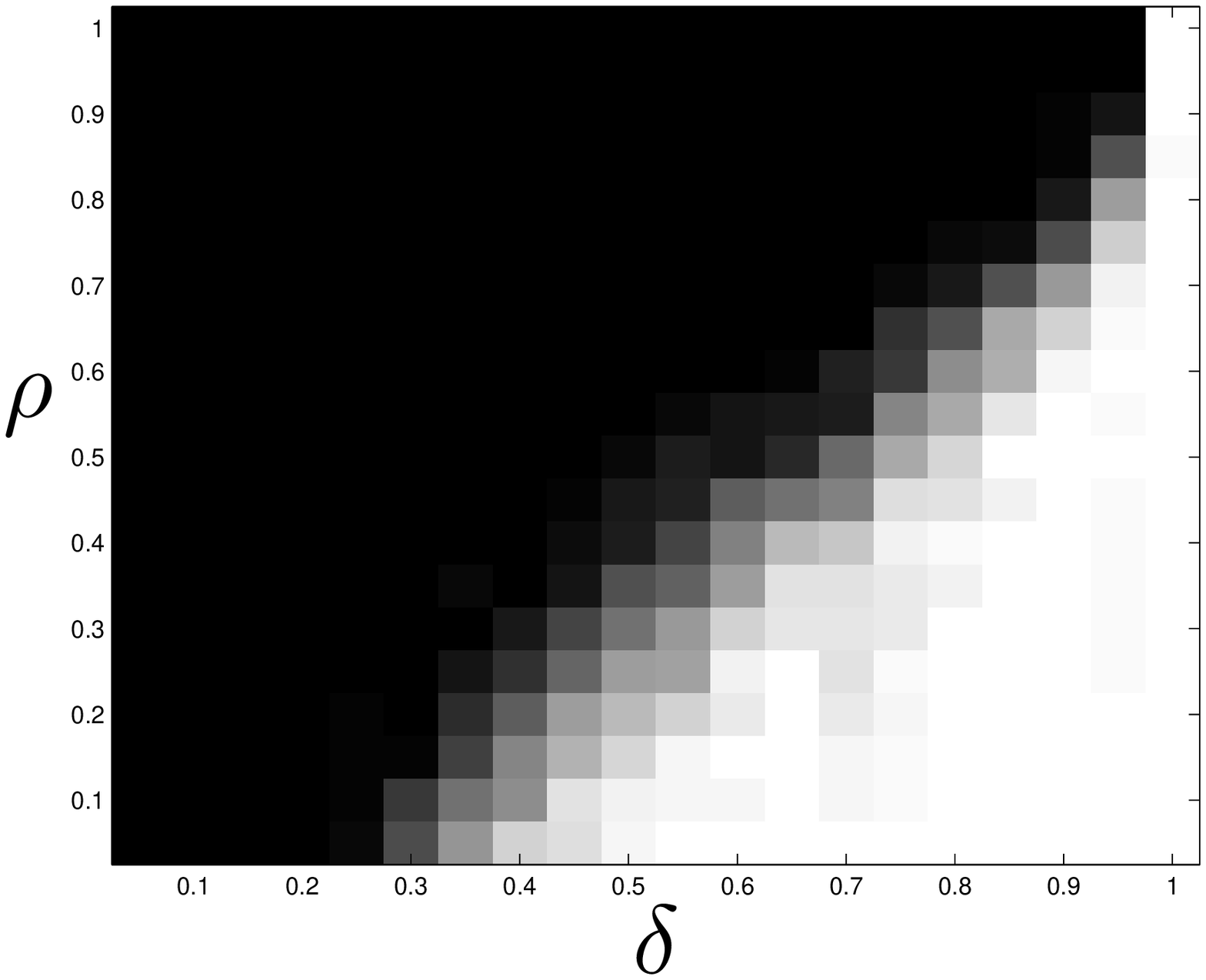}}%
\caption{Recovery rate for a random tight frame with $p=144$ and $d=120$. From left to right, up to bottom:
AIHT with an adaptive changing step-size,
AHTP with an adaptive changing step-size,
ACoSaMP,
ASP,
TDIHT with an adaptive changing step-size,
A-$\ell_1$-minimization and GAP.}
\label{fig:phaseDiagramAll1_2}
\end{figure*}

Figure~\ref{fig:phaseDiagramAll1_2} compares the reconstruction results of TDIHT with an adaptive changing step-size with the other algorithms.
It should be observed that TDIHT outperforms AIHT and provides comparable results with AHTP, where TDIHT is more computationally efficient than both of them. Though its performance is inferior to the other algorithms, it should be noted that its running time is faster by one order of magnitude compared to all the other techniques.
Therefore, in the cases that it succeeds, it should be favored over the other programs due to its better running time.

\begin{figure}[!t]
\centering
{\subfigure[Shepp-Logan phantom]{\includegraphics[width=1.7in]{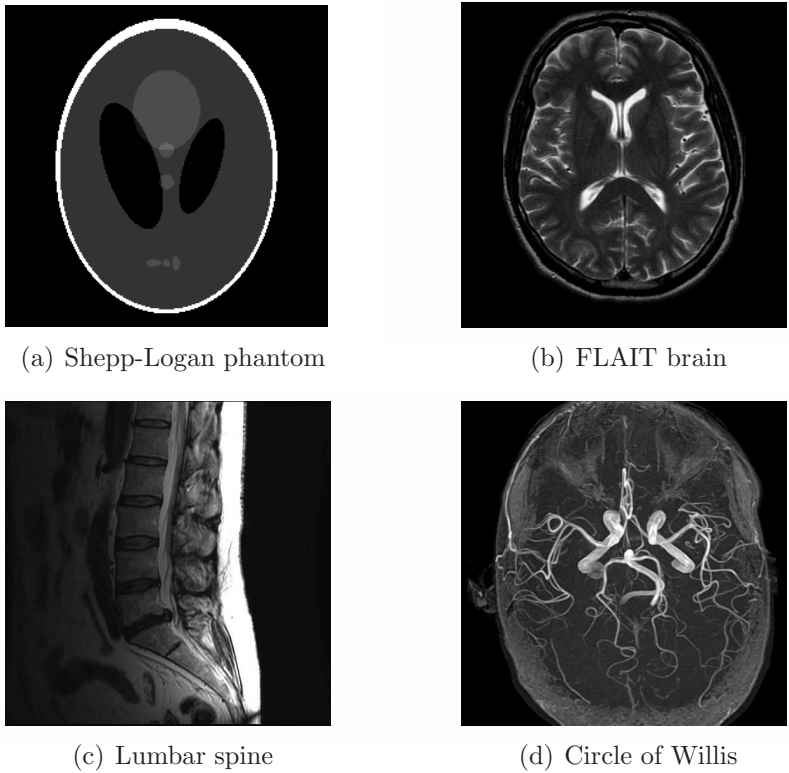} \label{fig:shepploganphantom}}%
\hfil
\subfigure[FLAIT brain]{\includegraphics[width=1.7in]{brain} \label{fig:brain}}%
\\ \subfigure[Lumbar spine]{\includegraphics[width=1.7in]{lumbspine} \label{fig:lumbspine}}%
\hfil
\subfigure[Circle of Willis]{\includegraphics[width=1.7in]{circle} \label{fig:circle}}}%
\caption{\rg{Test Images.}}
\label{fig:test_images}
\end{figure}

We turn now to test TDIHT for high dimensional signals.
\rg{We test the performance of several MRI images:  the \emph{Shepp-Logan phantom}, \emph{FLAIT brain} image, T2 Sagittal view of the \emph{lumbar spine} and the \emph{circle of Willis}. 
The first image is of size $256 \times 256$, while the other are of size $512 \times 512$. They are all presented in Fig.~\ref{fig:test_images}}.

We focus on the recovery of these images from a few number of Fourier measurements. 
With $\OM$ set to be the undecimated Haar transform with one level of resolution (redundancy four) and $\D$ its inverse transform,
we succeed to recover the phantom image using only $18$ sampled radial lines, which is only $6.5\%$ of the measurements.
This number is only slightly larger than the number needed for GAP, relaxed ASP (RASP) and Relaxed ACoSaMP (RACoSaMP) in \cite{Nam12Cosparse,Giryes14Greedy}. The advantage of TDIHT over these methods is its low complexity as it requires applying only $\M$ and its conjugate and $\OM$ and its inverse transform while in the other algorithms a high dimensional least squares minimization problem should be solved.
Note also that for AIHT and RAHTP the number of radial lines needed for recovery is $35$ and for IHT (with the decimated Haar operator with one level of resolution) we need more than $50$ radial lines.

\begin{figure}[!t]
\centering
{\subfigure[Naive PSNR = $18.3dB$]{\includegraphics[width=1.7in]{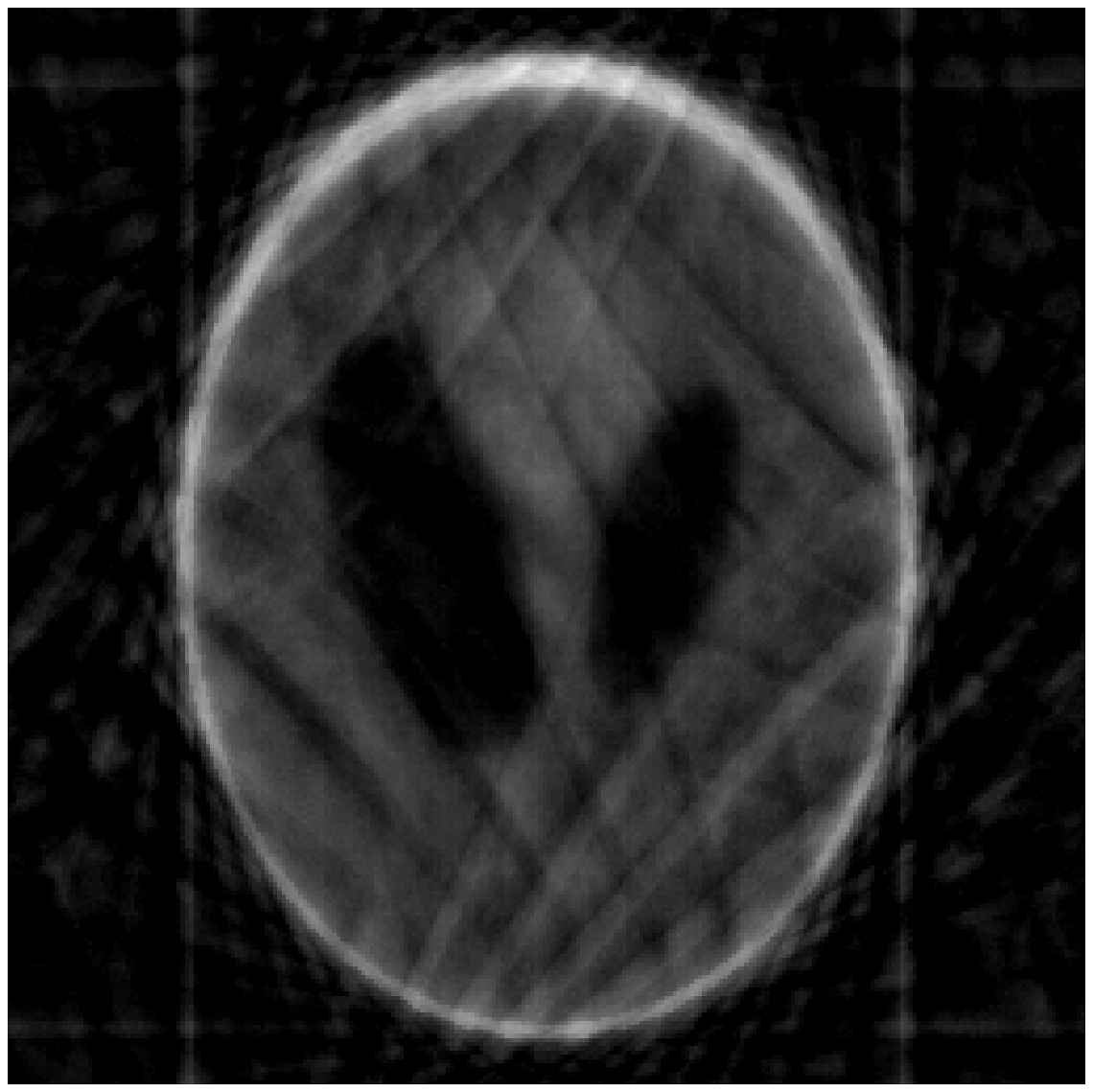} \label{fig:shepploganphantom_noisy}}%
\hfil
\subfigure[TDIHT PSNR = $36dB$]{\includegraphics[width=1.7in]{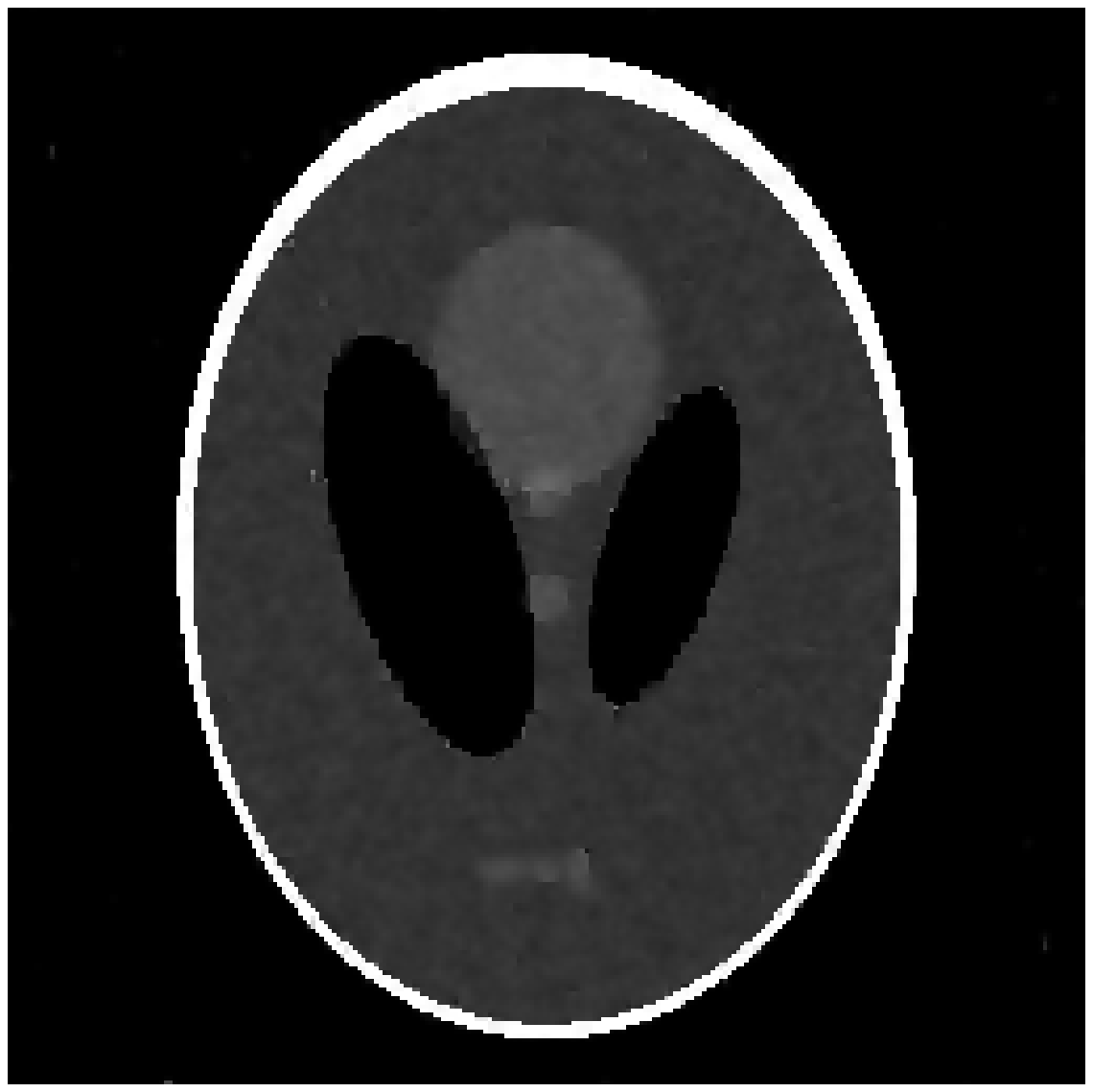} \label{fig:shepploganphantom_noisy_TDIHT}}}%
\caption{\rg{Recovery of the Shepp-Logan phantom from noisy $28$ radial lines with SNR of $20$. From left to right: Naive and TDIHT recovery.
Note that for the noiseless case TDIHT gets a perfect reconstruction using only $18$ radial lines.}}
\label{fig:shepploganphantom_reconstruction}
\end{figure}

Exploring the noisy case, we perform a reconstruction using TDIHT of a noisy measurement of the phantom with $28$ radial lines and signal to noise ratio (SNR) of $20$.
Figure~\ref{fig:shepploganphantom_noisy} presents the naive recovery from noisy image, the result of applying inverse Fourier transform on the measurements with zero-padding,
and Fig.~\ref{fig:shepploganphantom_noisy_TDIHT} presents the TDIHT reconstruction result.
We get a peak SNR (PSNR) of $36dB$. Note that RASP and GAPN require $22$ radial lines to get the same PSNR. However, this comes at the cost of higher computational complexity.

\begin{figure}[!t]
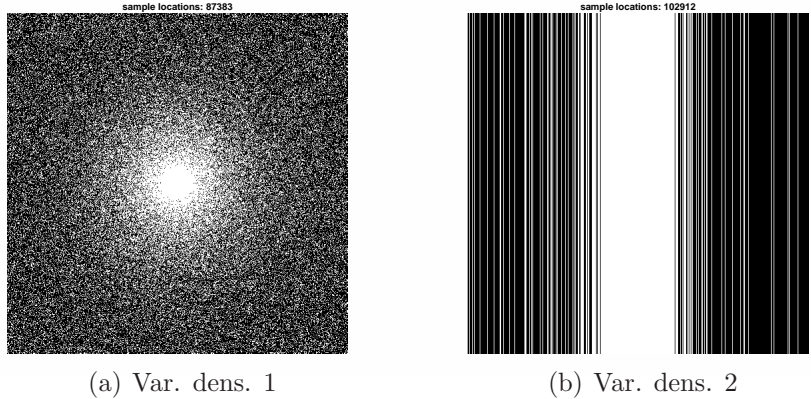

\centering
{\subfigure[Var. dens. 1]{\includegraphics[width=0.33\linewidth]{vardens} \label{fig:vardens1}}%
\hfil
\subfigure[Var. dens. 2]{\includegraphics[width=0.33\linewidth]{vertical} \label{fig:vertical}}}%
\caption{\rg{Pseudo-random variable-density undersampling patterns.}}
\label{fig:vardens}
\end{figure}


\begin{table}
\begin{center}
\begin{tabular}{|c|c|c|c|c|c|c|}
\hline
\bf Method &  \multicolumn{2}{c|}{\bf FLAIT Brain} & \multicolumn{2}{c|}{\bf Lumber Spine} & \multicolumn{2}{c|}{\bf The Circle of Willis} \\
\hline
 & Noiseless & Noisy & Noiseless & Noisy & Noiseless & Noisy  \\ \hline
Naive & 31.7 & 31.3 & 34.5 & 33.3  & 30 & 29.6 \\
TDIHT & 43.9 & 39.2 & 42.2 & 36.6 & 38.2 & 34.8 \\
RASP & 41.3 & 36.2  & 39.9 & 36.5  & 34.3 & 33.2\\
GAP & 42.1 & 35.8  & 40.4 & 34.6  & 36.2 &  32.7 \\
\hline
\end{tabular}
\end{center}
\caption{\rg{Recovery of FLAIT brain from var. dens. 1 samples, lumber spine from var. dens. 2 samples and circle of Willis from var. dens. 1 samples in the noiseless and noisy (with SNR 20) cases. The reconstruction quality is measured by PSNR.}}
\label{tbl:img_reconstruction}
\end{table}

\begin{figure}
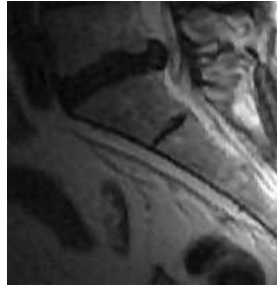
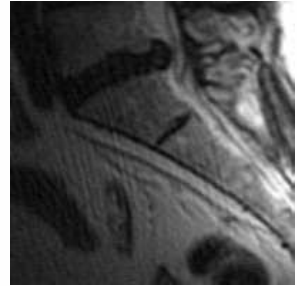
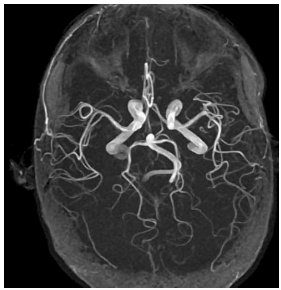
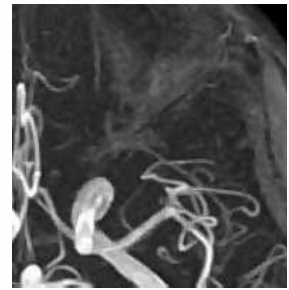
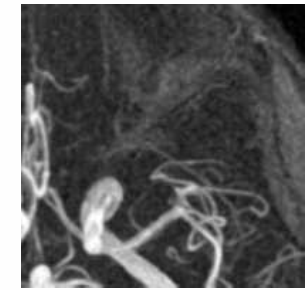

\centering
{\subfigure[TDIHT PSNR = $39.2dB$]{\includegraphics[width=0.275\linewidth]{brain_noisy_TDIHT}}%
\hfill
\subfigure[TDIHT zoom in]{\includegraphics[width=0.275\linewidth]{brain_noisy_TDIHT_part}}%
\hfill
\subfigure[Naive zoom in. PSNR = $31.3dB$]{\includegraphics[width=0.275\linewidth]{brain_noisy_naive_part} }%
\\
\subfigure[TDIHT PSNR = $36.6dB$]{\includegraphics[width=0.275\linewidth]{lumbspine_noisy_TDIHT} }
\hfill
\subfigure[TDIHT zoom in]{\includegraphics[width=0.275\linewidth]{lumbspine_noisy_TDIHT_part} }
\hfill
\subfigure[Naive zoom in. PSNR = $33.3dB$]{\includegraphics[width=0.275\linewidth]{lumbspine_noisy_naive_part} }
\\
\subfigure[TDIHT PSNR = $34.8dB$]{\includegraphics[width=0.275\linewidth]{circle_noisy_TDIHT} }
\hfill
\subfigure[TDIHT zoom in]{\includegraphics[width=0.275\linewidth]{circle_noisy_TDIHT_part} }
\hfill
\subfigure[Naive zoom in. PSNR = $29.6dB$]{\includegraphics[width=0.275\linewidth]{circle_noisy_naive_part}}}%
\caption{\rg{Reconstruction from noisy measurements with SNR 20. From left to right: TDIHT recovery, zoom-in on TDIHT recovery and zoom-in on naive recovery. From top to bottom: Recovery of FLAIT brain from var. dens. 1 samples, lumber spine from var. dens. 2 samples and circle of Willis from var. dens. 1. samples.}}
\label{fig:noisy_reconstruction}
\end{figure}

\rg{We perform similar experiments for the other images. Instead of uniformly sampling radial lines, we use pseudo-random variable-density undersampling patterns \cite{Lustig07Sparse}}\footnote{\rg{Unlike \cite{Lustig07Sparse} we perform our experiments with real (non-complex) images. The sampling patterns we use can be downloaded from \emph{http://www.eecs.berkeley.edu/$\sim${mlustig/CS.html}}.}} \rg{ presented in Fig.~\ref{fig:vardens}. We use var. dens. 1 with FLAIT brain and the circle of Willis and var. dens. 2 with lumber spine.}

\rg{As the images at hand are only approximately cosparse, we set a threshold, $0.01$ in the noiseless case and $0.05$ in the noisy one, such that each element in the cosparse representation below it is considered as zero. These thresholds create a model error in the recovery but provide a larger cosparsity value that eases the recovery. Moreover, in the noisy case, it is natural to set such a threshold since anyway small representation coefficients are being covered by the noise.}

\rg{Table~\ref{tbl:img_reconstruction} summarizes the recovery performance, in terms of PSNR, for the FLAIT brain, lumbar spine and circle of Willis images both for the noiseless and noisy cases.
To evaluate the performance of TDIHT in the noiseless case, we compare its PSNR with the one of the model error (which we get by applying $\OM$ on the original image followed by thresholding and multiplication by $\matr{D}$). For FLAIT brain and lumber spine we get errors, which are comparable, and even better for the latter, to their model errors $44.5dB$ and $44dB$. Note that such high PSNRs are equivalent to mean squared errors of the order of $10^{-5}$ and therefore we may say that we achieve a perfect recovery in these cases. 
This is not the case for the circle of Willis, in which the error is worse than the model error $41.3dB$. Note though that for this image, as well as for the other two images, we get better recovery error than the more sophisticated methods RASP and GAP.
Even in the noisy case, the quality we get with TDIHT is better than the one we have using those methods. }

\rg{Figs.~\ref{fig:noisy_reconstruction} presents the reconstruction outcome of TDIHT in the noisy case. To illustrate better the recovery gain, we present a zoom in on a part of the image and compare it to the naive recovery. 
The improvement can be seen clearly in all the three images. }

\section{Conclusion}
\label{sec:conc}

This paper presents a new algorithm, the transform domain IHT (TDIHT), for the cosparse analysis framework. In contrast to previous algorithms, TDIHT can be applied efficiently in high dimensional problems due to the fact that it requires applying only the measurement matrix $\matr{M}$, its transpose $\matr{M}^*$, the operator $\OM$ and its inverse $\matr{D}$ together with point-wise operations. The proposed algorithm is shown to provide stable recovery given that $\matr{M}$ has the D-RIP and that $\OM$ is a frame.
One of the limits of the proposed algorithm is that it assumes that the inverse of $\OM$ is given and can be easily applied in high dimensional problems. Though this is true for many types of operators like the discrete Fourier transform, wavelets and the Gabor frames \cite{Candes11Compressed}, it is not always possible to apply the inverse efficiently and there are even examples for operators for which such an inverse does not even exist, e.g., the finite difference operator \cite{Needell13Stable}.
A future work should pursue an efficient algorithm adapted to this case.


\appendix
\section{Proof of Lemma~\ref{lem:D_RIP_non_spars_up}}
\label{sec:D_RIP_non_spars_up_proof}

{\em Lemma~\ref{lem:D_RIP_non_spars_up}.}
Suppose that $\M$ satisfies the upper inequality of the D-RIP, i.e.,
\begin{eqnarray}
\norm{\M\D\w}_2 \le \sqrt{1+\delta_k}\norm{\D\w}_2 && \forall \w, \norm{\w}_0 \le k,
\end{eqnarray}
and that $\norm{\matr{D}}_2 \le \frac{1}{A}$.
Then for any representation $\w$ we have
\begin{eqnarray}
\norm{\M\D\w}_2 \le \frac{\sqrt{1+\delta_k}}{A}\left(\norm{\w}_2 + \frac{1}{\sqrt{k}}\norm{\w}_1 \right).
\end{eqnarray}

{\em Proof:}
We follow the proof of Proposition~3.5 in \cite{Needell09CoSaMP}.
We define the following two convex bodies
\begin{eqnarray}
&& \hspace{-0.3in} S = \conv\left\{\D\w : \w_{T^C} = 0, \abs{T} \le k, \norm{\D\w}_2 \le 1  \right\}, \\ &&
\hspace{-0.3in} K = \left\{\D\w : \frac{1}{A}\left(\norm{\w}_2 + \frac{1}{\sqrt{k}}\norm{\w}_1\right) \le 1 \right\}.
\end{eqnarray}
Since
\begin{eqnarray}
\norm{\M}_{S \rightarrow 2}  = \max_{\v \in S} \norm{\M\v}_2 \le
\sqrt{1+\delta_k},
\end{eqnarray}
it is sufficient to show that $\norm{\M}_{K \rightarrow 2} = \max_{\v \in K} \norm{\M\v}_2 \le  \norm{\M}_{S \rightarrow 2}$ which holds if $K \subset S$.
For proving the latter, let $\D\w \in K$ and $\{T_0, T_1, \dots, T_J\}$ be a set of distinct sets such that  $T_0$ is composed of the indexes of the  $k$-largest entries in $\w$, $T_1$ of the next $k$-largest entries, and so on.
Thus, we can rewrite $\D\w = \sum_{i=0}^J{\D\w_{T_i}} = \sum_{i=0}^J{\lambda_i\D\tilde{\w}_{T_i}}$, where $\lambda_i = \norm{\D\w_{T_i}}_2$ and $\tilde{\w}_{T_i} = \w_{T_i}/\lambda_i$.   Notice that by definition $\D\tilde{\w}_{T_i} \in S$. It remains to show that $\sum_{i=0}^J{\lambda_i} \le 1$ in order to show that $\D\w \in S$.
It is easy to show that $\norm{\w_{T_i}}_2 \le \frac{\norm{\w_{T_{i-1}}}_1}{\sqrt{k}}$.
Combining this with the fact that $\norm{\D}_2\le \frac{1}{A}$ leads to
\begin{eqnarray}
\sum_{i=1}^J{\lambda_i} & = & \sum_{i=1}^J{\norm{\D\w_{T_i}}_2} \le  \frac{1}{A}\sum_{i=1}^J{\norm{\w_{T_i}}_2} \\ \nonumber &\le & \frac{1}{A\sqrt{k}}\sum_{i=0}^{J-1}{\norm{\w_{T_i}}_1} \le \frac{1}{A\sqrt{k}}\norm{\w}_1.
\end{eqnarray}
Using the fact that $\lambda_0 = \norm{\D\w_{T_0}}_2 \le  \frac{1}{A}\norm{\w_{T_0}}_2 \le \frac{1}{A}\norm{\w}_2$, we have
\begin{eqnarray}
\sum_{i=0}^J{\lambda_i} \le \frac{1}{A}\left(\norm{\w}_2 + \frac{1}{\sqrt{k}}\norm{\w}_1 \right) \le 1,
\end{eqnarray}
where the last inequality is due to the fact that $\D\w \in K$.
\hfill $\Box$
\bigskip

\appendix
\section{Proof of Lemma~\ref{lem:TDIHT_xg_bound}}
\label{sec:TDIHT_xg_bound_proof}

{\em Lemma~\ref{lem:TDIHT_xg_bound}.}
Consider the same setup of Theorem~\ref{thm:TDIHT_frame_adversarial_noise}.
Then the $t$-th iteration of TDIHT satisfies
\begin{eqnarray}
&& \hspace{-0.3in} \norm{\left(\OM\x - \w_g^t\right)_{T \cup \hat{T}^t }}_2 \le 2\delta_{ak}\frac{B}{A}\norm{\left(\OM\x - \w_g^{t-1}\right)_{T \cup \hat{T}^{t-1} }}_2
\\ \nonumber&& \hspace{-0.15in}  +\frac{{1+\delta_{2k}}}{A}\left(\left(1+A + \frac{A^2}{2}\right)\norm{\OM_{T^C}\x}_2
+ \frac{1}{\sqrt{k}}\norm{\OM_{T^C}\x}_1 \right)
\\ \nonumber&&  \hspace{-0.15in}  + \norm{\OM_{T \cup  \hat{T}^t}\matr{M}^*\vect{e}}_2.
\end{eqnarray}

{\em Proof:}
Our proof technique is based on the one of IHT in  \cite{foucart10Sparse}, utilizing the properties of $\OM$ and $\D$.
Denoting $\w = \OM\x$ and using the fact that $\w_g^t = \OM\D\hat{\w}^{t-1} - \OM\M^* (\y-\M\D\hat{\w}^{t-1})$ we have
\begin{eqnarray}
&& \hspace{-0.4in}\norm{\left(\w - \w_g^t\right)_{T \cup \hat{T}^t }}_2 = \\ \nonumber && \hspace{-0.4in}
\norm{\left(\w - \OM\D\hat\w^{t-1}\right)_{T \cup \hat{T}^t } - \left(\OM\M^*(\y-\M\D\hat\w^{t-1})\right)_{T \cup \hat{T}^t }}_2.
\end{eqnarray}
By definition $\w = \OM\D\w$ and $\y = \M\D\w + \e$.
Henceforth
\begin{eqnarray}
\label{eq:z_zt_ineq1}
&& \hspace{-0.4in} \norm{\left(\w - \w_g^t\right)_{T \cup \hat{T}^t }}_2 \\ \nonumber && \hspace{-0.35in} =
\norm{\OM_{T \cup \hat{T}^t }\D\left(\w - \hat\w^{t-1}\right) - \OM_{T \cup \hat{T}^t }\M^*(\y-\M\D\hat\w^{t-1})}_2
\\ \nonumber && \hspace{-0.35in} = \norm{\OM_{T \cup \hat{T}^t }(\I - \M^*\M)\D\left(\w - \hat\w^{t-1}\right) - \OM_{T \cup \hat{T}^t }\M^*\e}_2 \\ \nonumber && \hspace{-0.35in} \le
\norm{\OM_{T \cup \hat{T}^t }(\I - \M^*\M)\D\left(\w_T - \hat\w^{t-1}\right) }_2 \\ \nonumber &&
+
\norm{\OM_{T \cup \hat{T}^t }(\I - \M^*\M)\D\w_{T^C} }_2
+
\norm{\OM_{T \cup \hat{T}^t }\M^*\e}_2,
\end{eqnarray}
where the last step is due to the triangle inequality.
Denote by $\P$ the projection onto $\range([\OM^*_{T\cup\hat{T}^{t}}, \D_{T \cup \hat{T}^{t-1}}])$, which is a subspace of vectors with $4k$-sparse representations. As $\w_T - \hat{\w}^{t-1}$ is supported on $T\cup T^{t-1}$ and $\M$ satisfies
the D-RIP for $[\OM^*, \D]$, we have
using norm inequalities and Lemma~\ref{lem:D_RIP_norm} that
\begin{eqnarray}
\label{eq:OM_IMM_D_z}
&& \hspace{-0.3in} \norm{\OM_{T \cup \hat{T}^t }(\I - \M^*\M)\D\left(\w_T - \hat\w^{t-1}\right) }_2
\\ \nonumber && \hspace{-0.1in} = \norm{\OM_{T \cup \hat{T}^t }\P(\I - \M^*\M)\P\D\left(\w_T - \hat\w^{t-1}\right) }_2
\\ \nonumber && \hspace{-0.1in} \le  \norm{\OM_{T \cup \hat{T}^t }}_2\norm{\P(\I - \M^*\M)\P }_2\norm{\D}_2\norm{\w_T - \hat\w^{t-1}}_2
\\ \nonumber && \hspace{-0.1in}
\le \delta_{4k}\frac{B}{A}\norm{\w_T - \hat\w^{t-1}}_2,
\end{eqnarray}
where $A$ and $B$ are the frame constants and we use the fact that $\norm{\D}_2\le \frac{1}{\norm{\OM}_2} \le \frac{1}{A}$ and that $\norm{\OM_{T\cup \hat{T}^{t}}}_2 \le \norm{\OM}_2 \le B$.
Notice that when $\OM$ is tight frame, $\OM^* = \D$ and thus  $\OM_T^* = \D_T$. Hence, we have $\delta_{3k}$ instead of $\delta_{4k}$ since
 $\range([\OM^*_{T\cup\hat{T}^{t}}, \D_{T \cup \hat{T}^{t-1}}])$ is a subspace of vectors with $3k$-sparse representations.

For completing the proof we first notice that $\w_T - \hat{\w}^{t-1} = (\w - \hat{\w}^{t-1})_{T \cup \hat{T}^{t-1}}$ and that $\hat{\w}^{t-1}$ is the best $k$-term approximation of $\w_g^{t-1}$ in the $\ell_2$ norm sense. In particular it is also the best $k$-term approximation of
$(\w_g^{t-1})_{T\cup \hat{T}^{t-1}}$ and therefore $\norm{(\hat\w^{t-1} - \w_g^{t-1} )_{T\cup \hat{T}^{t-1}}}_2 \le \norm{(\w - \w_g^{t-1})_{T\cup \hat{T}^{t-1}}}_2$. Starting with the triangle inequality and then applying this fact we have
\begin{eqnarray}
\label{eq:z_zt1_ineq_z_zg}
&& \hspace{-0.3in} \norm{\w_T - \hat\w^{t-1}}_2
\\ \nonumber &&  = \norm{(\w - \w_g^{t-1} + \w_g^{t-1} - \hat\w^{t-1})_{T\cup \hat{T}^{t-1}}}_2
\\ \nonumber &&\le  \norm{(\w - \w_g^{t-1})_{T\cup \hat{T}^{t-1}}}_2 +\norm{(\hat\w^{t-1} - \w_g^{t-1} )_{T\cup \hat{T}^{t-1}}}_2 \\ \nonumber &&\le  2\norm{(\w - \w_g^{t-1})_{T\cup \hat{T}^{t-1}}}_2.
\end{eqnarray}
Combining \eqref{eq:z_zt1_ineq_z_zg} and \eqref{eq:OM_IMM_D_z} with \eqref{eq:z_zt_ineq1}
leads to
\begin{eqnarray}
\label{eq:z_zgt_ineq_semifin}
\hspace{-0.15in}\norm{\left(\w - \w_g^t\right)_{T \cup \hat{T}^t }}_2 &\le &
 2\delta_{4k}\frac{B}{A}\norm{(\w - \w_g^{t-1})_{T\cup \hat{T}^{t-1}}}_2
\\ \nonumber && +
\norm{\OM_{T \cup \hat{T}^t }(\I - \M^*\M)\D\w_{T^C} }_2
\\ \nonumber && +
\norm{\OM_{T \cup \hat{T}^t }\M^*\e}_2.
\end{eqnarray}
It remains to bound the second term of the rhs. By using the triangle inequality
and then the D-RIP with the fact that $\norm{\OM_{T \cup \hat{T}^t }\matr{D}}_2 \le \norm{\OM\matr{D}}_2 \le 1$ we have
\begin{eqnarray}
\label{eq:OM_IMM_DzTC_ineq}
&& \hspace{-0.4in} \norm{\OM_{T \cup \hat{T}^t }(\I - \M^*\M)\D\w_{T^C} }_2
\\ \nonumber && \le
\norm{\OM_{T \cup \hat{T}^t }\D\w_{T^C} }_2  +\norm{\OM_{T \cup \hat{T}^t }\M^*\M\D\w_{T^C} }_2
\\ \nonumber &&\le
\norm{\w_{T^C} }_2  +\sqrt{1+\delta_{2k}}\norm{\M\D\w_{T^C} }_2
\\ \nonumber && \le
\norm{\w_{T^C} }_2  +\frac{{1+\delta_{2k}}}{A}\left(\norm{\w_{T^C}}_2 + \frac{1}{\sqrt{k}}\norm{\w_{T^C}}_1 \right),
\end{eqnarray}
where the last inequality is due to Lemmas~\ref{lem:k_inequality} and \ref{lem:D_RIP_non_spars_up}.
The desired result is achieved by plugging \eqref{eq:OM_IMM_DzTC_ineq} into \eqref{eq:z_zgt_ineq_semifin}
and using the fact that $1+\frac{1+\delta_{2k}}{A} \le \frac{(1+\delta_{2k})(1+A)}{A}$.
\hfill $\Box$
\bigskip

\section*{Acknowledgment}
The author would like to thank Michael Elad, Yoram Bresler and Yaniv Plan for fruitful discussions.
R. Giryes is grateful to the Azrieli Foundation for the award of
an Azrieli Fellowship. 
Flait Brain, lumber spine and circle of Willis images are taken from {http://www3.americanradiology.com}.
The authors would like to thank the anonymous reviewers for their helpful and constructive comments that greatly contributed to improving the final version of the paper.


\bibliographystyle{elsarticle-num}
\bibliography{../analysis}

\end{document}